\documentclass[final]{siamltex}

\usepackage{color}
\usepackage{textcomp}
\usepackage[mathcal]{euscript}
\usepackage{mathrsfs}
\usepackage{amsfonts,amssymb,amsbsy, amsmath}
\usepackage{graphicx,empheq}
\usepackage{nomencl}
\usepackage[page,header]{appendix}
\usepackage{url}
\makeglossary
\newcommand{\g}{{ \bf g}}
\newcommand{\y}{{ \bf y}}

\newcommand{\z}{{\bf z}}
\newcommand{\cLtwo}{{\mathcal L}_2 (0,L)}

\usepackage{pifont}
\usepackage[tableposition=top]{caption}

\newcommand{\mc}{\mathcal}
\newcommand{\mb}{\mathbf}

\def\e1{{\varepsilon_{1}}}
\def\b1{{\beta_{11}}}
\def\bp3{{\beta_{33}}}
\def\ep3{{\varepsilon_{3}}}

\newtheorem{thm}{Theorem}[section]

\newtheorem{Lem}{Lemma}[section]

\MHInternalSyntaxOn
\providecommand*\phantomword[3][c]{%
\mathchoice
{\MT_phantom_word:NNnn #1\displaystyle {#2}{#3}}%
{\MT_phantom_word:NNnn #1\textstyle {#2}{#3}}%
{\MT_phantom_word:NNnn #1\scriptstyle {#2}{#3}}%
{\MT_phantom_word:NNnn #1\scriptscriptstyle {#2}{#3}}%
}
\def\MT_phantom_word:NNnn #1#2#3#4{%
\@begin@tempboxa\hbox{$\m@th#2#4$}%
\setlength\@tempdima{\widthof{$\m@th#2#3$}}%
\hbox{\hb@xt@\@tempdima{\csname bm@#1\endcsname}}%
\@end@tempboxa}
\MHInternalSyntaxOff

\title{Modeling and related results for current-actuated piezoelectric beams by including magnetic effects\footnote{T\lowercase{he financial support of  the }NSERC  D\lowercase{iscovery }G\lowercase{rant Program for this research is gratefully acknowledged.} }
}

\author{Kirsten A. Morris  \thanks{Department of Applied Mathematics, University of Waterloo, Waterloo, ON N2L3G1, Canada ({\tt kmorris@uwaterloo.ca}).}
        \and Ahmet \"{O}zkan \"{O}zer \thanks{Department of Mathematics, University of Nevada, Reno, NV 89503, USA({\tt aozer@unr.edu}).}}

\begin{document}

\maketitle
\begin{abstract}
Piezo-electric material can be controlled with current as the electrical variable, instead of voltage. The main purpose of this paper is to derive the governing equations for a current-controlled piezo-electric beam and to investigate stabilizability.
Besides the consideration of current control, there are several new aspects to the model here.
Most significantly, magnetic effects are included.
 For the electromagnetic part of the model,  electrical potential and magnetic vector potential  are chosen to be quadratic-through thickness to include the induced effects of the electromagnetic field.  Two sets of decoupled system of partial differential equations are obtained; one for stretching motion and another one for bending motion.  Hamilton's principle is used to derive a boundary value problem that models a single piezo-electric beam actuated by a charge (or current) source  at the electrodes.  Current or charge controllers at the electrodes can only control the stretching motion. Attention is therefore focused on  control of the stretching equations in this paper. It is shown that the Lagrangian of the beam is invariant under certain transformations. A Coulomb-type gauge condition which is widely used in the electromagnetic theory is used here. This gauge condition decouples the electrical potential equation from the equations of the magnetic potential. A semigroup approach is used to prove that the Cauchy problem is well-posed. Unlike the voltage or charge actuation, a bounded control operator in the natural energy space is obtained in the current actuation case.
 The paper concludes with analysis of stabilizability and comparison with other actuation approaches and models.
 \end{abstract}

 \section{Introduction}
\begin{figure}[h]
\centering
\includegraphics[width=3.6in]{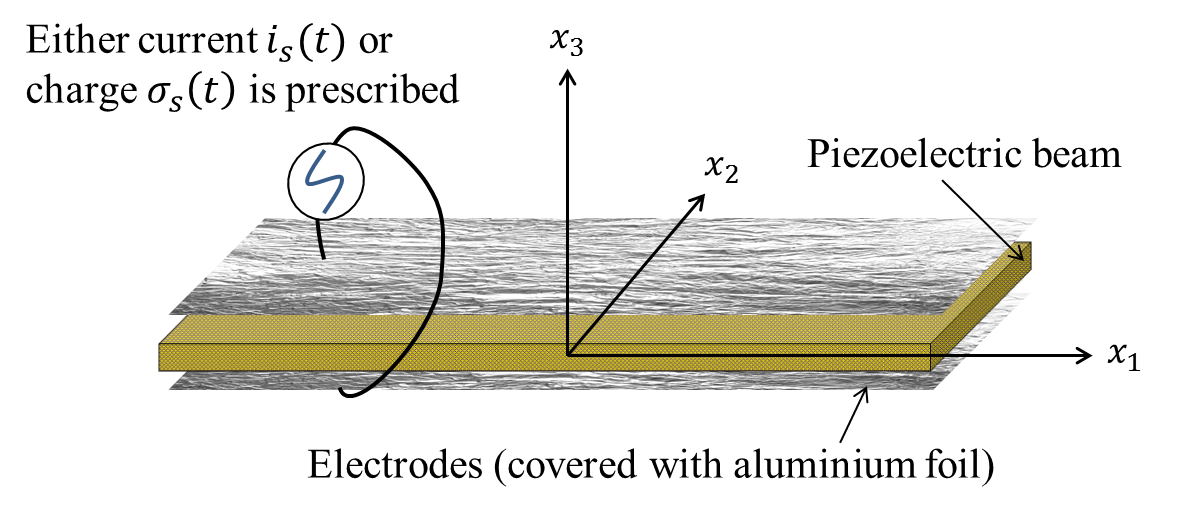}
\caption{\footnotesize When either charge $\sigma_s(t)$ or current $i_s(x,t)$ is prescribed to the electrodes, an electric field is created between the electrodes, and therefore the beam/plate either shrinks or extends. Unlike the voltage actuation, the input-output hysteresis reduces substantially.}
\label{pbeam}
\end{figure}

 Piezoelectric materials are elastic beam/plates covered by electrodes at the top and bottom surfaces, insulated at the edges (to prevent fringing effects), and connected to an external electric circuit. (See Figure \ref{pbeam}.) They convert mechanical energy to electrical and \emph{magnetic energy}, and vice versa.    These materials are widely used in civil, aeronautic and space structures due to their small size and high power density. These materials can be actuated by either external mechanical forces or electrical forces.  There are mainly three ways to (electrically) actuate piezoelectric materials: voltage, current or charge.  Piezoelectric materials have been traditionally activated by a voltage source \cite{Cao-Chen,Chee,Dest,Rogacheva,Smith,Stanway,Tiersten,Tzou}. It is well-known that the control operator is unbounded in the energy space if the piezoelectric structure is controlled by a voltage or a charge source, for instance see \cite{Banks-Smith,Dest,Hansen,O-M,Smith}.

  Hysteresis  occurs  in the voltage-strain relationship for piezo-electric structures; see for instance, \cite{Smith}. This complicates control of these materials.   One way to avoid  hysteresis is by applying only low voltages, but this prevents these structures from being used at their maximum potential. Therefore controller design needs to consider hysteresis in order to obtain maximum accuracy and effectiveness. Some approaches are passivity \cite{rob} and inverse compensation \cite{SmithII}.  Another way to reduce the hysteresis is current or charge actuation, see for instance \cite{Comstock,F-M,Hagood,MGN,MG1,Newcomb-Flinn}. Existing models for current control use only  circuit equations attached to the standard elastic beam  equations, and magnetic effects are not considered.  Magnetic effects, even though  small, were previously shown to be very important to the controllability and stabilizability of piezoelectric beams
\cite{Yang1}. In fact, it has been  shown \cite{O-M} that voltage-controlled beams cannot even be asymptotically stabilized for certain material parameters.

 In this paper,  dynamic magnetic effects are included in the derivation of a model for  piezoelectric beams actuated by a current (or charge) source.  The electromagnetic field is described in terms of  scalar electric potential and magnetic vector potential, After deriving expressions the various contributions to the energy, Hamilton's Principle is used to derive
 a system of partial differential equations modelling  the coupling between the mechanical  and the electro-magnetic  dynamics.
 These equations do not have a unique solution since the potentials are not uniquely determined. This is because the  Lagrangian corresponding to Maxwell's equations  is invariant under certain transformations; for instance see \cite[page 80]{L-D-I}).
Obtaining a system of equations with a unique solution requires a appropriate gauge condition. A number of  gauges are possible.
A Coulomb type of transformation is used here.  Implementation of this transformation  simplifies the equations considerably.  The original highly coupled system of equations becomes a system of equations where the equations corresponding to the electrical variables are completely decoupled from the ones involving magnetic potential variables.
Well-posedness of the model in an appropriate Hilbert space  is then established. The norm in the Hilbert space corresponds to the energy of the system. It is shown that the spectrum of the generator consists entirely of imaginary eigenvalues. Stabilizability of the model is  compared to voltage control, as well as to the case where magnetic effects are neglected.
Some of the results presented in this paper, in particular,
Lemma 2.1, Theorem 3.1, Theorems 4.1, 4.2 and a weaker version of Theorem \ref{main-thm2}
were previously reported in the conference paper  \cite{O-M3}.

\begin{table}[h]

 {\footnotesize{
\renewcommand{\arraystretch}{1.3}
   \begin{tabular}[c]{|c|l|c|l|}
\hline
$A$ & Magnetic potential vector & $\rho$ & Mass density per unit volume \\
  \hline
  $B$ &  Magnetic flux density vector  &   $n$ & Surface unit outward normal vector \\
    \hline
  $c, \alpha$ &  Elastic stiffness coefficients  &  $\sigma_s$ & Surface charge density \\
    \hline
$D$ &  Electric displacement vector& $\sigma_b$ & Volume charge density \\
    \hline
  $E$ & Electric field intensity vector &  $S$ & Strain tensor  \\
    \hline
   $\varepsilon$ &  Permittivity coefficients      & $T$ & Stress tensor \\
    \hline
$h$ & Thickness of the beam     &  $U_1$ & $x_1$ component of the displacement field    \\
    \hline
$i_b$ & Volume current density & $U_3$ & $x_3$ component of the displacement field\\
  \hline
$i_s$ & Surface  current density    &   $v$ & Longitudinal disp. of the centerline of the beam \\
  \hline
 $\phi$ & Electric potential   & $V$ & Voltage \\
    \hline
$\gamma$ &  Piezoelectric coefficients  & $w$ & Transverse displacement of the beam\\
    \hline
 $\mu$ & Magnetic permeability of beams  &  & \\
    \hline
\end{tabular}
}}

\caption{Notation}
\end{table}

\section{Physical Model}
\label{Sec-II}
Let $x_1, x_2$ be the longitudinal directions, and $x_3$ be the transverse directions (see Figure \ref{pbeam}). Assume that the piezoelectric beam occupy the region $\Omega=[0,L]\times [-r,r]\times [-\frac{h}{2}, \frac{h}{2}]$ with the boundary $\partial \Omega,$   the electroded region and the insulated region, where $L> > h.$
 Throughout this paper,  dots denote differentiation with respect to time, that is $\dot{x}(t)=\frac{d x}{dt} .$

  A very widely-used linear constitutive relationship \cite{Tiersten} for piezoelectric beams  is
\begin{eqnarray}
\label{cons-eqq10}
\left( \begin{array}{l}
 T \\
 D \\
 \end{array} \right)=
\left[ {\begin{array}{*{20}c}
   c & -\gamma^{\text{T}}  \\
   \gamma & \varepsilon  \\
\end{array}} \right]\left( \begin{array}{l}
 S \\
 E \\
 \end{array} \right)
\end{eqnarray}
where $T=(T_{11}, T_{22}, T_{33}, T_{23}, T_{13}, T_{12})^{\rm T}$ is the stress vector, \\$S=(S_{11}, S_{22}, S_{33}, S_{23}, S_{13}, S_{12})^{\rm T}$  is the strain vector, $D=(D_1, D_2, D_3)^{\text T}$ and  $E=(E_1, E_2, E_3)^{\text{T}}$ are the electric displacement  and the electric field vectors, respectively, and moreover, the matrices $[c], [\gamma], [\varepsilon]$ are the matrices with elastic, electro-mechanic and dielectric constant entries (for more details see \cite{Tiersten}). Under the assumption of transverse isotropy and polarization in $x_3-$direction, these matrices reduce to

\begin{eqnarray}
\nonumber &c= \left( {\begin{array}{*{20}c}
   {c_{11} } & {c_{12} } & {c_{13} } & 0 & 0 & 0  \\
   {c_{21} } & {c_{22} } & {c_{23} } & 0 & 0 & 0  \\
   {c_{31} } & {c_{32} } & {c_{33} } & 0 & 0 & 0  \\
   0 & 0 & 0 &    c_{44} & 0 & 0 \\
      0 & 0 & 0 &    0 & c_{55} & 0 \\
         0 & 0 & 0 &    0 & 0 & c_{66}
\end{array}} \right),~~~~  \Gamma= \left( {\begin{array}{*{20}c}
   0 & 0 & 0 &    0 & \gamma_{15} & 0 \\
      0 & 0 & 0 &    -\gamma_{15} & 0 & 0 \\
         \gamma_{31} & \gamma_{31} & \gamma_{33} &    0 & 0 & 0
\end{array}} \right)&\\
\nonumber  &\varepsilon= \left( {\begin{array}{*{20}c}
   \varepsilon_{11} & 0 & 0 \\
      0 & \varepsilon_{22} & 0\\
         0 &  0 & \varepsilon_{33} \\
\end{array}} \right)&
\end{eqnarray}
We assume that all forces acting in the $x_2$ direction are zero which implies a beam. Moreover, $ T_{33}$ is also assumed to be zero. Therefore
$$T=(T_{11}, T_{13})^{\text T}, S=(S_{11}, S_{13})^{\text T}, D=(D_1, D_3)^{\text T}, E=(E_1, E_3)^{\text{T}}$$ and (\ref{cons-eqq10}) reduces to
\begin{eqnarray}
\nonumber
\left( \begin{array}{l}
  T_{11} \\
  T_{13} \\
 D_1 \\
 D_3
 \end{array} \right)
&= \left( {\begin{array}{*{20}c}
   {c_{11} }  &  0 & 0 & -\gamma_{31}  \\
   0 & {c_{55} } & -\gamma_{15}  & 0 \\
   0 & \gamma_{15} & \varepsilon_{11}&   0 &  \\
      \gamma_{31} & 0 & 0  & \varepsilon_{55}  \\
\end{array}} \right)\left( \begin{array}{l}
  S_{11} \\
  S_{13} \\
 E_1 \\
 E_3
 \end{array} \right).
\end{eqnarray}
Let $(U_1, U_3)$  denote the displacement field vector of a point $(x_1, x_3).$  Continuing with the Euler-Bernoulli beam small-displacement assumptions, the displacement field is given as the following
\begin{eqnarray}
\label{kirc} && U_1=v-x_3\frac{\partial w}{\partial x_1}, \quad  U_3=w
\end{eqnarray}
where $v=v(x_1)$ and $w=w(x_1)$ denote the longitudinal displacement of the center line in $x_1$ direction, and transverse    displacement of the beam, respectively.  Since $S_{13}=\frac{1}{2}\left(\frac{\partial U_1}{\partial x_3} + \frac{\partial U_3}{\partial x_1}\right)=0,$ the only nonzero strain component is given by
\begin{eqnarray}S_{11}=\frac{\partial U_1}{\partial x_1}=\frac{\partial v}{\partial x_1}-x_3 \frac{\partial^2 w}{\partial x_1^2}.\label{strains}
\end{eqnarray}
To keep the notation simple let \begin{eqnarray}
\alpha=c_{11}, \quad \gamma=\gamma_{31}, \quad \gamma_1=\gamma_{15},\quad \varepsilon_1=\varepsilon_{11},\quad \varepsilon_3=\varepsilon_{33}.
\end{eqnarray}
With the new notation,  the linear constitutive equations for an Euler-Bernoulli piezoelectric beam are
\begin{eqnarray}
 \label{cons-eq10} &&\left\{
  \begin{array}{ll}
     T_{11}=\alpha S_{11}-\gamma E_3& \\
  T_{13} = -\gamma_{1} E_1& \\
   D_1=\varepsilon_{1} E_1 &\\
  D_3=\gamma S_{11}+\varepsilon_{3}E_3&
  \end{array} \right.
\end{eqnarray}

 Let $\mb K, \mb P, \mb E$ and  $\mb B$  be kinetic, potential, electrical, and magnetic energies of the beam, respectively, and let $\mb W$ be the work done by the external forces. To model charge or current-actuated piezoelectric beams we use the following Lagrangian \cite{Lee}
\begin{eqnarray}\label{Lag}  \mb{ L}= \int_0^T \left[\mb{K}-(\mb{P}-\mb{E}+\mb B)+ \mb{W}\right]dt\end{eqnarray}
for which we use the constitutive equations (\ref{cons-eq10}) where the pair $(S,E)$ belongs to the set of independent variables.  In the above, $\mb P-\mb E+\mb B$ is called electrical enthalpy. Note that in modeling piezoelectric beams by voltage-actuated electrodes we use a different Lagrangian  \begin{eqnarray}\label{tildeL}  \tilde{\mb L}= \int_0^T \left[\mb{K}-(\mb{P}+\mb{E})+\mb B +\mb {W}\right]~dt\end{eqnarray} for which the constitutive equations (\ref{cons-eq10}) are written in terms of the independent variables $(S,D).$ The Lagrangian  $\tilde{\mb L}$ can be obtained by applying a Legendre transformation to $\mb L.$ Here $\mb P+\mb E$ denotes the total stored energy of the beam, and $\mb B$ acts as the electrical kinetic energy of the beam. This case is studied in \cite{O-M}. Depending on the prescribed quantity at the electrodes, Lagrangian is chosen either $\mb L$ or $\tilde{\mb L}.$

 The full set of Maxwell's equations is; see for example  \cite[Page 332]{Duvaut-L}),
  \begin{eqnarray}
\nonumber &\nabla\cdot D =~\sigma_b  \quad{\rm{in}} \quad\Omega \times \mathbb{R}^+~\quad&\text{(Electric Gauss's ~law)}\\
\nonumber &\nabla\cdot B=~0 \quad{\rm{in}} \quad\Omega \times \mathbb{R}^+~\quad  & \text{(Gauss's law of magnetism)}\\
\nonumber &\nabla\times E=~-\dot B  \quad{\rm{in}} \quad\Omega \times \mathbb{R}^+~\quad&\text{(Faraday's law)}\\
\label{Maxwell}  &\frac{1}{\mu}(\nabla\times B)= ~i_b +  \dot D  \quad{\rm{in}} \quad\Omega \times \mathbb{R}^+~\quad& \text{(Amp\'{e}re-Maxwell law)}
\end{eqnarray}
with one of the essential electric boundary conditions prescribed on the electrodes
\begin{eqnarray}
\label{charge_boun}  - D\cdot n =~\sigma_s(t)  & \quad{\rm{on}} \quad {\partial\Omega} \times \mathbb{R}^+~\quad&\text{(Charge )}\\
\label{current_boun} \frac{1}{\mu}(B \times n) =~i_s(t) & \quad{\rm{on}} \quad {\partial\Omega} \times \mathbb{R}^+~\quad & \text{(Current)}\\
\label{voltage_boun} \phi=~V(t)  & \quad{\rm{on}} \quad {\partial\Omega} \times \mathbb{R}^+~\quad& \text{(Voltage)}
\end{eqnarray}
and appropriate mechanical boundary conditions at the edges of the beam  (the beam is clamped, hinged, free, etc.).
Here $B$ denotes the magnetic field vector, and $\sigma_b, i_b, \sigma_s, \\i_s, V, \mu, n$ denote body charge density, body current density, surface charge density,
 surface current density, voltage, magnetic permeability, and unit normal vector to the surface $\partial\Omega,$ respectively. In this paper we consider  only current and charge-driven electrodes (i.e. we ignore(\ref{voltage_boun})). The voltage-driven electrode case is handled in details in \cite{O-M}. In modeling piezoelectric beams, there are mainly three approaches including electric and magnetic effects \cite{Tiersten}:
\vspace{0.1in}

\begin{description}
  \item[i) Electrostatic electric field:] Electrostatic electric field approach is the most widely-used approach in the literature.
  It completely ignores magnetic effects: $B=\dot D=i_b=\sigma_b=0.$  Maxwell's equations  (\ref{Maxwell}) reduce to
   $\nabla\cdot D =0$ and $\nabla \times E=0.$ Therefore,  there exist a scalar electric potential
   such that $E=-\nabla \phi$  and $\phi$ is determined up to a constant.

  \item[ii) Quasi-static electric field:] This approach rules out some of the magnetic effects (non-magnetizable materials) \cite{Tiersten}:
  $\sigma_b=i_b=0.$ However, $\dot D$ and $B$ are non-zero. Therefore, (\ref{Maxwell}) reduce to
      $$ \nabla \cdot D=0, ~~~\nabla\cdot B=0, ~~~\dot B= -\nabla\times E, ~~~ \dot D= \frac{1}{\mu}(\nabla\times B).$$
      The equation $\nabla\cdot B=0$
       implies that there exists a vector $A$ such that $B=\nabla\times A.$ This vector is called the {\em magnetic potential .}
        It follows from substituting $B$ to $\dot B= -\nabla\times E$ that there exists a scalar electric potential $\phi$ such that \begin{equation}
\label{imp1}E=-\dot A-\nabla\phi.
\end{equation}
 where $\dot A$ stands for the induced electric field due to the time-varying magnetic effects.
 One simplification in this approach is to ignore $A$ and $\dot A$ since $A, \dot A\ll \phi.$  With this assumption  $\dot D$ may be non-zero.

  \item[iii) Fully dynamic electric field: ] Unlike the quasi-static assumption,  $A$ and $\dot A$ are left in the model.  Depending on the type of material, body charge density $\sigma_b$ and body current density $i_b$ can also be  non-zero. Note that even though the  displacement current $\dot D$ is assumed to be non-zero in both quasi-static and fully dynamic approaches, the term $\ddot D$  is  zero in quasi-static approach since $\dot A=0.$
\end{description}
Since the piezoelectric materials are not perfectly  insulated, the electric field $E$ causes currents to flow when conductivity occurs. Therefore the time-dependent equation of the continuity
of electric charge must be employed. In this paper, we follow the fully dynamic approach  to include all of the magnetic effects. If we take the divergence of both sides of Amp\'{e}re-Maxwell equation (\ref{Maxwell}), we obtain
\begin{eqnarray}\frac{1}{\mu} \nabla\cdot \left(\nabla\times B\right)= \nabla \cdot i_b +  \nabla \cdot \dot D.\label{cont}\end{eqnarray}
The term on the left hand side of the equation above  is zero, and therefore by using Gauss's law (\ref{Maxwell}), we obtain the following electric continuity condition
\begin{eqnarray}\dot \sigma_b+ \nabla \cdot i_b =0 \quad {\text in} ~~\Omega
\label{elec-cont1}\end{eqnarray}
The physical interpretation of (\ref{elec-cont1}) is the local conservation of electrical charge. From (\ref{Maxwell})
$$\frac{1}{\mu} \int_{\partial\Omega}\left(\nabla\times B\right)\cdot n~ dS  = \int_{\partial\Omega} \left(i_b \cdot n +   \dot D \cdot n\right)~dS, $$
and use the charge boundary conditions (\ref{charge_boun}) with $i_s(x,t)\equiv 0$
\begin{eqnarray}\nonumber 0= \frac{1}{\mu} \int_{\partial\Omega} \nabla \times B \cdot n~dS &=&  \int_{\partial\Omega}\left( i_b^3 -  \dot\sigma_s\right)~dS
 \end{eqnarray}
 or, alternatively,  use the current boundary condition (\ref{current_boun}) with $\sigma_s\equiv 0$
 \begin{eqnarray}\nonumber \frac{1}{\mu} \int_{\partial\Omega} \nabla \times B \cdot n~dS &=&  \int_{\partial\Omega} i_b^3 ~dS
 \end{eqnarray}
 where $n$ is the outward unit normal vector on $\partial\Omega.$
Hence we obtain surface electric continuity conditions (or compatibility conditions)
\begin{eqnarray}\dot \sigma_s - i_b^3= 0, \quad {\rm or,}\quad\frac{ di_s}{dx}- i_b^3= 0 \quad \text{on}~ \partial\Omega.
\label{elec-cont2}\end{eqnarray}
For more details, the reader can refer to \cite[Section 3.9]{Eringen}.

 Henceforth, to simplify the notation, $x=x_1$ and $z=x_3.$


\label{Fully-dy} Note that piezoelectricity is the direct result of piezoelectric effect,  which is, once the external electric field is applied to the electrodes, strain is produced and therefore the beam/plate extends or shrinks (direct effect), whereas, when the plate/beam extends and shrink, it produces electric voltage which is so-called the induced (inverse) effect. The linear through-thickness assumption of the electric potential $\phi(x,z)=\phi^0(x)+z\phi^1(x)$ completely ignores the induced potential effect since $\phi$ is completely known as a function of voltage. For example, when the voltage is prescribed at the electrodes, i.e.  $\phi\left(\frac{h}{2}\right)=V$ and $\phi\left(-\frac{h}{2}\right)=0,$  we have
 $$\phi=\frac{V}{2}+z\frac{V}{h},$$ and therefore the electric field component in the transverse direction $E_3$ becomes uniform in the transverse direction, i.e. $E_3=-\phi^1=-\frac{V}{h},$ as we consider electrostatic and quasi-static assumptions. Therefore the induced effect is ignored in the linear-through thickness assumption. In this paper, we use a quadratic-through thickness potential distribution that takes care of the induced effect and improves the modeling accuracy:
 \begin{eqnarray}
 \label{scalarpot} && \phi(x,z)=\phi^0(x)+z \phi^1(x)+ \frac{z^2}{2} \phi^2(x).
 \end{eqnarray}

Since we are in the beam theory, we assume that the magnetic vector potential $A$ has nonzero components only in $x$ and $z$ directions. To keep the consistency with $\phi$, we assume that $A$ is quadratic through-thickness as well:
\begin{eqnarray}
 \label{vectorpot}&& A(x,z) =
\left( \begin{array}{c}
 A_1(x,z) \\
 0 \\
 A_3(x,z) \\
 \end{array} \right)
=\left( \begin{array}{c}
 A_1^0(x) + z A_1^1(x)  + \frac{z^2}{2} A_1^2(x)\\
 0 \\
  A_3^0(x) + z A_3^1(x)  + \frac{z^2}{2} A_3^2(x)\\
 \end{array} \right).
\end{eqnarray}
By (\ref{imp1})
\begin{eqnarray}
\nonumber &&E_1= - \left( \dot A_1^0 +z \dot A_1^1 + \frac{z^2}{2} \dot A_1^2 \right)-\left(({\phi}^0)_x+z (\phi^1)_x+ \frac{z^2}{2} (\phi^2)_x\right),\\
\label{e3} &&E_3= -\left( \dot A_3^0 + z \dot A_3^1 +  \frac{z^2}{2} \dot A_3^2\right)-\left( \phi^1 + z \phi^2 \right).
\end{eqnarray}
Now we use the constitutive equations (\ref{cons-eq10}) along with (\ref{kirc}), (\ref{strains}), and (\ref{scalarpot})-(\ref{e3}) to write 
\begin{eqnarray}
\nonumber \mb{E}-\mb{P}&=&\frac{1}{2}\int_\Omega \left(D_1 E_1+D_3 E_3- T_{11}S_{11}-T_{13}S_{13}\right)~dX \\
\nonumber &=&\frac{1}{2}\int_\Omega \left( -\alpha S_{11}^2 +2\gamma S_{11}E_3 + \e1 E_1^2 + \ep3 E_3^2\right)~dX \\
\nonumber &=& \frac{1}{2} \int_0^L \left[-{\alpha} h\left(v_x^2+ \frac{h^2}{12}w_{xx}^2\right)- 2{\gamma}h  \left(\left(\phi^1+ \dot A_3^0 + \frac{h^2}{24} \dot A_3^2\right)v_x-\frac{h^2}{12}w_{xx}\left(\phi^2+\dot A_3^1\right)\right) \right.\\
\nonumber && \quad  + {\e1}h \left( (\phi^0_x)^2+\frac{h^2}{12} (\phi^1_x)^2+ \frac{h^4}{320} (\phi^2_x)^2 +(\dot A_1^0)^2+ \frac{h^2}{12} (\dot A_1^1)^2+ \frac{h^4}{320} (\dot A_1^3)^2  \right)\\
\nonumber && \quad  + {\ep3}h \left( (\phi^1)^2+\frac{h^2}{12} (\phi^2)^2 +(\dot A_3^0)^2+ \frac{h^2}{12} (\dot A_3^1)^2+ \frac{h^4}{320} (\dot A_3^2)^2  \right)\\
\nonumber && \quad + 2{\e1}h \left( (\phi^0)_x \dot A_1^0 + \frac{h^2}{24} (\phi^0)_x( \phi^2)_x + \frac{h^2}{24} \dot A_1^0 \dot A_1^2+ \frac{h^2}{24} (\phi^0)_x \dot A_1^2 \right.\\
\nonumber &&\left. \quad\quad\quad\quad\quad+ \frac{h^2}{24} (\phi^2)_x \dot A_1^0 +   \frac{h^2}{12} (\phi^1)_x \dot A_1^1+  \frac{h^4}{320} (\phi^2)_x \dot A_1^2  \right)\\
\label{E-P} && \quad\left. \quad  + 2{\ep3}h \left( \phi^1 \dot A_3^0 + \frac{h^2}{24} \dot A_3^0 \dot A_3^2 + \frac{h^2}{24} \phi^1  \dot A_3^2 + \frac{h^2}{12} \phi^2 \dot A_3^0 \right)\right]~dx,
\end{eqnarray}
\begin{eqnarray}
\nonumber \mb{B}&=& \frac{\mu}{2}\int_\Omega (\nabla \times A) \cdot  (\nabla \times A) ~dX \\
\nonumber &=& \frac{\mu}{2} \int_0^L\int_{-h/2}^{h/2} \left(A_1^1 + z A_1^2  -(A_3^0)_x-z(A_3^1)_x-\frac{z^2}{2}(A_3^2)_x\right)^2 ~dzdx \\
\nonumber &=& \frac{\mu h}{2} \int_0^L\left[ (A_1^1)^2+ \frac{h^2}{12} (A_1^2)^2 +((A_3^0)_x)^2+\frac{h^2}{12}((A_3^1)_x)^2 + \frac{h^4}{320}((A_3^2)_x)^2 \right. \\
\label{BB} && \left.-2\left( A_1^1 ~(A_3^0)_x +\frac{h^2}{24} A_1^1(A_3^2)_x - \frac{h^2}{12} A_1^2(A_3^1)_x -\frac{h^2}{24}(A_3^0)_x(A_3^2)_x\right) \right]~dx,\quad\quad\quad
\end{eqnarray}
\begin{eqnarray}
\label{KK} \mb{K}&=& \frac{\rho}{2} \int_\Omega \left(\dot U_1^2+ \dot U_3^2\right)~dX= \frac{\rho h}{2} \int_0^L \left(\dot v^2 + \dot w^2 +\frac{h^2}{12}\dot w_x^2\right)~dx,
\end{eqnarray}

Now we define the work $\mb W$ done by the external forces. We first define the body force resultants  $i_b, \sigma_b, i_s, \sigma_s$ as in  \cite{Lagnese-Lions}:
 \begin{eqnarray}\nonumber &&i_b= \int_{-h/2}^{h/2} \tilde i_b~ dz , \quad \sigma_b= \int_{-h/2}^{h/2} \tilde \sigma_b ~ dz  \\
  \nonumber &&i_s= \int_{-h/2}^{h/2} \tilde i_s~ dz , \quad \sigma_s = \int_{-h/2}^{h/2} \tilde \sigma_s~ dz.\end{eqnarray}
In the above the surface charge density $\tilde \sigma_s$ and surface current density $\tilde i_s$ are  independent of $z$ since they are given at the electrodes.  For the Euler-Bernoulli beam, it is appropriate to assume that body charge $\tilde \sigma_b$ and body current  $\tilde i_b$ are independent of $z,$
 $$i_b=\tilde i_b h,\quad \sigma_b = \tilde \sigma_b h, \quad i_s = \tilde i_s h,\quad \sigma_s=\tilde \sigma_s h.$$
 We choose either surface charge $\sigma_s $ or $i_s$ to be non-zero, and  either $i_b$ or $\sigma_b$ to be nonzero, depending on the type of actuation.  We assume that there are no mechanical external forces acting on the beam.  The work done by the electrical external forces is \cite{Lee}
\begin{eqnarray}
\nonumber \mb{W}&=& \int_{\Omega} \left(~-\tilde \sigma_b~ \phi+\tilde i_b \cdot A  \right) ~ dX + \int_{\partial \Omega} \left(-\tilde\sigma_s~\phi + \tilde i_s\cdot  A\right)   ~ dX\\
\nonumber &=&\int_{\Omega} \left(~-\tilde \sigma_b ~\phi+ \tilde i_b^1 A_1 \right) ~ dX + \int_{\partial \Omega} \left(-\tilde \sigma_s~\phi + \tilde  i_s^1\cdot  A_1\right)   ~ dX\\
\nonumber  &=& -\int_0^L\int_{-h/2}^{h/2} \tilde \sigma_b \left(\phi^0(x)+z \phi^1(x)+ \frac{z^2}{2} \phi^2(x)\right) ~dzdx\\
 \nonumber && + \int_0^L\int_{-h/2}^{h/2} \tilde i_b^1 \left(A_1^0(x) + z A_1^1(x)  + \frac{z^2}{2} A_1^2\right) ~ dz dx \\
 \nonumber && + \int_0^L  \left(-\tilde \sigma_s \left(\phi(h/2)-\phi(-h/2)\right) + \tilde i_s^1\left(A_1(h/2)-A_1(-h/2)\right)\right)  ~dx \\
\label{work} &=& \int_0^L  \left(  -\sigma_b\left(\phi^0 + \frac{h^2}{24}\phi^2\right)-\sigma_s+i_b^1\left(A_1^0 + \frac{h^2}{24}A_1^2\right)- \sigma_s ~\phi^1 +i_s^1 ~A_1^1 \right)~ dx
\end{eqnarray}
where  $i_s(x,t)=(i_s^1(x,t),0,0),$ and $i_b(x,t)=(i_b^1(x,t), 0, 0). $  In the above $i_s$ has only one nonzero component since  $i_s \perp B $, and $i_s \perp n$ by (\ref{current_boun}). Moreover, $i_b$ has only one nonzero component since we assumed that there is no force acting in the $x_2$ and $x_3$ directions.

 If the magnetic effects are neglected, a variational approach cannot be used in the case of current actuation since $A\equiv 0$ and so $\mb W\equiv 0.$
This is very different from the charge and voltage actuation cases since for charge and voltage actuation $\mb W$ is not a function of $A.$

\section{Derivation of Governing Equations}
Assume that both ends of the piezoelectric beam are free. The application of Hamilton's principle, setting the variation of Lagrangian $\mb L$ in  (\ref{Lag}) with respect to the all kinematically  admissible displacements $$\{v, w,\phi^0, \phi^1, \phi^2, A_1^0, A_1^1, A_1^2, A_3^0, A_3^1, A_3^2\}$$ to zero, yields stretching equations and bending equations  respectively:
\begin{equation}\label{Stretch-EE}\left\{
\begin{array}{l l}
    \rho h  \ddot v -{\alpha}h v_{xx} -{\gamma}h \left((\phi^1)_x +  (\dot A_3^0)_x  + \frac{h^2}{24} (\dot A_3^2)_x \right)= 0 & \\
    -\frac{{\e1} h^3}{12}\left((\phi^1)_{xx}+ (\dot A_1^1)_x\right) +{\ep3}h\left( \dot A_3^0  +\frac{h^2}{24} \dot A_3^2+ \phi^1\right)-{\gamma}hv_x= \sigma_s & \\
    \frac{{\e1}h^3}{12} \ddot A_1^1 + \frac{{\e1}h^3}{12} (\dot \phi^1)_x-\mu h\left((A_3^0)_x +\frac{h^2}{24}(A_3^2)_x- A_1^1\right)=i_s^1  & \\
    {\ep3}h \left( \ddot A_3^0 + \frac{h^2}{24} \ddot A_3^2  + \dot \phi^1\right)-\mu h\left((A_3^0)_{xx}+\frac{h^2}{24}(A_3^2)_{xx}-(A_1^1)_x\right)-{\gamma}h\dot v_x=0& \\
     \frac{{\ep3}h^3}{24} \left( \ddot A_3^0  + \dot \phi^1\right) + \frac{\ep3 h^5}{320} \ddot A_3^2-\mu h^3\left(\frac{h^2}{24}(A_3^0)_{xx}+\frac{h^2}{320}(A_3^2)_{xx}-\frac{1}{24}(A_1^1)_x\right)&\\
     \quad\quad\quad\quad\quad -\frac{{\gamma}h^3}{24}\dot v_x=0&
\end{array}\right.
\end{equation}
\begin{equation}\label{Bending-E}\left\{
\begin{array}{l l}
    \rho h  \ddot w - \frac{\rho h^3}{12} \ddot w_{xx} + \frac{{\alpha} h^3}{12} w_{xxxx} -\frac{{\gamma} h^3}{12} \left( (\phi^2)_{xx}) +  (\dot A_3^1)_{xx}  \right)= 0 & \\
   -{\e1}h \left( (\phi^0)_{xx}+\frac{h^2}{24}(\phi^2)_{xx}+(\dot A_1^0)_x + \frac{h^2}{24}(\dot A_1^2)_x\right)  = \sigma_b & \\
    -\frac{{\e1} h^3}{24} \left( (\phi^0)_{xx}+\frac{h^2}{24}(\phi^2)_{xx}+(\dot A_1^0)_x + \frac{h^2}{24}(\dot A_1^2)_x\right)+\frac{{\gamma}h^3}{24}w_{xx}& \\
    \quad\quad\quad\quad\quad-\frac{{\e1}h^5}{720}(\phi^2)_{xx}- \frac{{\e1}h^5}{720}(\dot A_1^2)_x  +\frac{{\ep3}h^3}{12}\left( \phi^2+ \dot A_3^1\right)=\frac{h^2 \sigma_b}{24} & \\
    {\e1}h \left( \ddot A_1^0  + \frac{h^2}{24} \ddot A_1^2 +(\dot\phi^0)_x + \frac{h^2}{24}(\dot\phi^2)_x \right) =i_b^1  &\\
   \frac{{\e1}h^3}{24} \left( \ddot A_1^0 + \frac{h^2}{24}(\ddot A_1^2)+(\dot\phi^0)_x + \frac{h^2}{24}(\dot\phi^2)_x\right) & \\
\quad\quad\quad\quad\quad + \frac{{\e1}h^5}{720} \ddot A_1^2 +\frac{{\e1}h^5}{720}(\dot \phi^2)_x+\frac{h^3\mu}{12}\left( A_1^2-(A_3^1)_x \right)=\frac{h^2 i_b^1}{24} &\\
\frac{{\ep3}h^3}{12} \left( \ddot A_3^1 - \frac{\mu}{{\ep3}}( A_3^1)_{xx} \right)+   \frac{{\ep3} h^3}{12} \dot\phi^2 +\frac{h^3\mu}{12}(A_1^2)_x+\frac{{\gamma}h^3}{12}\dot w_{xx}=0 &
\end{array}\right.
\end{equation}
with the natural boundary  conditions at $x=0,L$
\begin{small}
\begin{equation}\label{table:bcss} \left\{
\begin{array}{r r}
 {\alpha}hv_x +{\gamma}h\left(\phi^1+\dot A_3^0+\frac{h^2}{24}\dot A_3^2 \right)=0 & {\rm {(Lateral ~ force)}} \\
  \frac{h^3}{12}\left(-{\alpha}w_{xx} +{\gamma}\phi^2\right)=0 &{\rm{(Bending ~moment)}} \\
 -\rho\ddot w_x +{\alpha}w_{xxx} -{\gamma}(\phi^2)_x=0 & {\rm {(Shear) }} \\
   \e1 h\left( \dot A_1^0+\frac{h^2}{24} \dot A_1^2+({\phi}^0)_x+ \frac{h^2}{24} (\phi^2)_x\right)=0 & {\rm {(Charge) }}\\
   \frac{\e1 h^3}{12}\left(\dot A_1^1+ (\phi^1)_x\right)=0 & {\rm {(First~ charge ~   moment)}}\\
  {\e1}h^3\left(\frac{1}{12} \dot A_1^0+\frac{h^2}{160} \dot A_1^2 + \frac{1}{12}({\phi}^0)_x+ \frac{h^2}{160}(\phi^2)_x\right)=0 & {\rm { (Second ~charge~ moment) }} \\
   \mu h\left( A_1^1 -(A_3^0)_x-\frac{h^2}{24}(A_3^2)_x\right)=0 & {\rm {(Current) }}\\
 \frac{\mu h^3}{12}\left( A_1^2 - (A_3^1)_x\right)=0 & {\rm {(First ~current ~moment)}}  \\
   \mu h^3\left(\frac{1}{24} A_1^1  -\frac{1}{24}(A_3^0)_x-\frac{h^2}{320}(A_3^2)_x\right)=0 & {\rm {(Second ~current~ moment)}}
\end{array}\right.
\end{equation}
\end{small}
The bending motion is described by the Rayleigh beam equation coupled to the electromagnetic equations. If the rotational inertia of the cross section of the beam is ignored, the terms $\ddot w_{xx}$ in (\ref{Bending-E})and $\ddot w_x$ in (\ref{table:bcss}) go away.

The last equation in (\ref{Stretch-EE}) can be simplified by using the previous one to get
\begin{equation}\label{Stretch-E}\left\{
\begin{array}{l l}
    \rho h  \ddot v -{\alpha}h v_{xx} -{\gamma}h \left((\phi^1)_x +  (\dot A_3^0)_x  + \frac{h^2}{24} (\dot A_3^2)_x \right)= 0 & \\
    -\frac{{\e1} h^3}{12}\left((\phi^1)_{xx}+ (\dot A_1^1)_x\right) +{\ep3}h\left( \dot A_3^0  +\frac{h^2}{24} \dot A_3^2+ \phi^1\right)-{\gamma}hv_x= \sigma_s & \\
    \frac{{\e1}h^3}{12} \ddot A_1^1 + \frac{{\e1}h^3}{12} (\dot \phi^1)_x-\mu h\left((A_3^0)_x +\frac{h^2}{24}(A_3^2)_x- A_1^1\right)=i_s^1  & \\
    {\ep3}h \left( \ddot A_3^0 + \frac{h^2}{24} \ddot A_3^2  + \dot \phi^1\right)-\mu h\left((A_3^0)_{xx}+\frac{h^2}{24}(A_3^2)_{xx}-(A_1^1)_x\right)-{\gamma}h\dot v_x=0& \\
      \frac{{\ep3}h^3}{24} \left( \ddot A_3^0 + \frac{h^2}{24} \ddot A_3^2 + \dot \phi^1\right)-\frac{\mu h^3}{24}\left((A_3^0)_{xx}+\frac{h^2}{24}(A_3^2)_{xx}-(A_1^1)_x\right) &\\
    \quad\quad\quad\quad\quad -\frac{{\gamma}h^3}{24}\dot v_x+ \left(\frac{{\ep3} h^5}{720} \ddot A_3^2 -\frac{\mu h^5}{720}(A_3^2)_{xx}\right)=0
\end{array}\right.
\end{equation}
Note that the stretching  (\ref{Stretch-E}) and bending (\ref{Bending-E})  equations are completely decoupled when only one type of external electrical force is present.  It will be assumed, as is common in practice, that there is no free body charge or current,  i.e. $\sigma_b \equiv i_b\equiv 0. $ Then the bending equations (\ref{Bending-E}) are entirely uncontrolled and also decoupled from the  stretching equations (\ref{Stretch-E}). Therefore, from this point on,  only the stretching equations (\ref{Stretch-E}) are considered with the corresponding boundary conditions at $x=0,L$
\begin{equation}\label{a} \left\{
\begin{array}{l r}
 {\alpha}hv_x +{\gamma}h\left(\phi^1+\dot A_3^0+\frac{h^2}{24}\dot A_3^2 \right)=0 & {\rm {(Lateral ~ force)}} \\
   \frac{\e1 h^3}{12}\left(\dot A_1^1+ (\phi^1)_x\right)=0 & {\rm {(First~ charge ~   moment)}}\\
   \mu h\left( A_1^1 -(A_3^0)_x-\frac{h^2}{24}(A_3^2)_x\right)=0 & {\rm {(Current) }}\\
  \frac{ \mu h^3}{24}\left(A_1^1  -(A_3^0)_x-\frac{h^2}{24}(A_3^2)_x\right)-\frac{\mu h^5}{720}(A_3^2)_x=0 & {\rm {(Second ~current~ moment)}}
\end{array}\right.
\end{equation}
The last two boundary conditions can also be simplified as \begin{eqnarray}\label{aa}\left\{A_1^1 -(A_3^0)_x\right\}_{x=0,L}=\left\{(A_3^2)_x\right\}_{x=0,L}=0.
\end{eqnarray}


\label{Sec-III}
 The magnetic potential vector $A$ and the electric potential $\phi$ are not uniquely defined by (\ref{Maxwell}).
  In fact, the Lagrangian $\mb L$  (\ref{Lag}) is invariant under a large class of transformations.

 \begin{thm}\label{invariant}
 For any scalar $C^1$ function $\chi=\chi(x, z, t),$  the Lagrangian $\mb L$ is invariant under the  transformation
\begin{eqnarray}
\nonumber  & A \mapsto \tilde A:= A+\nabla \chi& \\
\label{Gaugemap} & \phi \mapsto \tilde \phi:= \phi-\dot\chi.&
\end{eqnarray}
\end{thm}
\textbf{Proof:} By (\ref{Gaugemap}),  $\tilde A$ and $\tilde \phi$ satisfy $$\tilde B=\nabla \times \tilde A= \nabla \times A+\nabla \times \nabla\chi=\nabla \times A=B$$  $$\tilde E=-\dot{\tilde A}-\nabla \tilde\phi=-\dot A-\nabla \dot \chi-\nabla\phi + \nabla \dot \chi=-\dot A- \nabla \phi=E.$$ This implies that $\mb E- \mb P$ and $\mb B$ defined respectively by (\ref{E-P}) and (\ref{BB}) are invariant under the transformation. Since $\mb K$ in (\ref{KK}) is independent of $A$ and $\phi,$ we need to check if $\mb W$ defined by (\ref{work}) is invariant under (\ref{Gaugemap}).  We choose the arbitrary scalar function $\chi$ to be quadratic-through thickness $\chi=\chi^0+ z\chi^1 +\frac{z^2}{2}\chi^2$ to be consistent with the choices of  $\varphi$ and $A$ in (\ref{scalarpot}) and (\ref{vectorpot}), respectively.  We also assume that $\chi$ satisfies the stationary conditions $\left|\chi^0=\chi^1=\chi^2\right|_{t=0,T}=0$ for compatibility. By (\ref{Gaugemap}) we have
\begin{eqnarray}
\nonumber &&\tilde A_1^0= A_1^0+(\chi^0)_x, \quad A_1^1= A_1^1+(\chi^1)_x,\quad A_1^2= A_1^2+(\chi^2)_x,\quad \tilde A_3^0= A_3^0 + \chi^1, \quad \tilde A_3^2=A_3^2, \\
\nonumber &&\tilde \phi^0=\phi^0 - \dot \chi^0, \quad \tilde\phi^1=\phi^1-\dot \chi^1, \quad \tilde\phi^2=\phi^2-\dot \chi^2,
\end{eqnarray}
and therefore
\begin{eqnarray}
\nonumber \int_0^T {\tilde {\mc{W}}}~dt &=& \int_0^T \int_0^L  \left(  -\sigma_b\left({\tilde {\phi^0}} + \frac{h^2}{24}\tilde \phi^2\right)-\sigma_s~\tilde{\phi^1}+i_b^1\left(\tilde A_1^0 + \frac{h^2}{24}\tilde A_1^2\right) +i_s^1 ~\tilde A_1^1  \right)~ dx dt\\
\nonumber &=&\int_0^T \mc W dt  + \int_0^T\int_0^L  \left( \sigma_b\left(\dot \chi^0 +\frac{h^2}{24}\dot\chi^2\right)+\sigma_s\dot \chi^1 \right.\\
\nonumber &&\quad \left.+ \left(i_b^1\left((\chi^0)_x + \frac{h^2}{24}(\chi^2)_x \right)  +i_s^1 ~(\chi^1)_x\right)\right) ~dxdt.\\
\nonumber &=&\int_0^T {{\mc{W}}}~dt +  \int_0^T \int_0^L \left[\left(\frac{di_s^1}{dx} \right) \chi^1 -\frac{di_b^1}{dx}\left(\chi^0 + \frac{h^2}{24}\chi^2 \right) \right]~dxdt \\
\nonumber && +\int_0^T \int_0^L \left(\dot \sigma_b\left( \chi^0 +\frac{h^2}{24}\chi^2\right)+\dot \sigma_s \chi^1 \right) ~dx dt\\
\nonumber && + h \left[\int_0^T \left(i_s^1~ \chi^1 +i_b^1 \left(\chi^0 + \frac{h^2}{24}\chi^2\right)\right) ~dt\right]_0^L \\
\nonumber &&+\left[\int_0^L \left(\sigma_b\left( \chi^0 +\frac{h^2}{24}\chi^2\right)+\sigma_s \chi^1 \right) ~dt\right]_0^T     \\
\nonumber &=&\int_0^T { {\mc{W}}}~dt
\end{eqnarray}
where we used $i_s^1=i_b^1=0$ at the insulated edges of electrodes. Hence, the Lagrangian $\mb L$ is  invariant under the transformation (\ref{Gaugemap}). $\square$

Since $\mb L$ is invariant under  transformations of type (\ref{Gaugemap}),  the electric and magnetic potentials are not uniquely determined by (\ref{Stretch-E}) and (\ref{a}). An additional condition can be added to remove the ambiguity. The additional condition is generally known as  a {\em gauge} and it is generally chosen to simplify the equations.
It is often convenient to choose the gauge to decouple the electrical potential equation  from the equations of the magnetic potential.
The Coulomb gauge in standard electromagnetic theory is defined by
 $$\rm{Div} A=0\quad {\rm in}\quad \Omega,\quad  A\cdot n=0\quad {\rm on}\quad \partial\Omega.$$
 This is one of the gauges commonly used in electromagnetic theory. With this additional condition, the Maxwell equations (\ref{Maxwell}) written in terms of the potentials
 \begin{eqnarray}&&-\nabla^2 \phi -\frac{\partial (\nabla \cdot A)}{\partial t}=0 \label{eqm1}\\
 && \frac{ \partial^2A}{\partial t^2}-\nabla^2 A=-\nabla\frac{\partial \phi}{\partial t}-\nabla (\nabla \cdot A) \nonumber
\end{eqnarray}
  are decoupled (for instance see \cite[page 80]{L-D-I}) and  (\ref{eqm1}) becomes
  $$ -\nabla^2 \phi = 0 .$$
  Here it was assumed that $A$ and $\phi$ are quadratic in the thickness variable $z$. So in  (\ref{E-P})-(\ref{work}) integration by parts is
 with respect to the $x$ variable, but not $z$.  Thus, (\ref{Stretch-E}) is not identical to (\ref{eqm1}).
 Examining equation (\ref{Stretch-E} b) the appropriate condition to decouple the magnetic and electric equations  is
\begin{eqnarray}
\label{Colombo}& -\frac{\e1 h^2}{12}(A_1^1)_x + \ep3\left( A_3^0+\frac{h^2}{24} A_3^2 \right)= 0 .
\end{eqnarray}
The boundary conditions
\begin{equation}\label{gauge-bc}(A_1^1)(0)=(A_1^1)(L)=0\end{equation}
are added so that the boundary conditions (\ref{a})-(\ref{aa}) are also decoupled.

The gauge condition uniquely determines $\phi$ and $A.$ Let $A$, $\tilde{A}$ be  potentials that satisfy (\ref{Colombo})  and (\ref{gauge-bc}), and also are related by a transformation  of the form (\ref{Gaugemap}). Then the arbitrary scalar function $\chi$ has to satisfy the following differential equation
\begin{eqnarray}\nonumber 0 &=& -\frac{\e1 h^2}{12}(\tilde A_1^1)_x + \ep3\left(\tilde A_3^0 +\frac{h^2}{24} \tilde A_3^2\right)\\&=&-\frac{\e1 h^2}{12}\left(( A_1^1)_{x} + (\chi^1)_{xx}\right)
 + \ep3\left(A_3^0+\chi^1 + \frac{h^2}{24} A_3^2\right)\\
\label{eq4} &=& -\frac{\e1 h^2}{12}(\chi^1)_{xx}+\ep3 \chi^1.
\end{eqnarray}
with the boundary conditions  $$(\chi^1)_x(0)=(\chi^1)_x(L)=0$$
where we used (\ref{gauge-bc}).
 Since (\ref{eq4}) with this boundary condition has only the trivial solution $\chi^1\equiv 0$ it follows that the additional conditions (\ref{Colombo},\ref{gauge-bc}) uniquely define the potentials $\phi$ and $A$ in (\ref{Gaugemap}). The existence and uniqueness of the solutions of the system (\ref{Stretch-E}) with  (\ref{Colombo},\ref{gauge-bc}) will be analyzed in details in Section \ref{Sec-IV}.

Note that other components of $\chi$ are coupled through the bending equation (\ref{Bending-E}). Showing that these components are equal to zero requires to choice of   another gauge condition similar to (\ref{Colombo}). Since  bending equations are not considered in this paper, this point is not considered.

Define $\eta:= A_3^0 + \frac{h^2}{24} A_3^2, \quad \theta:=A_1^1.$
The gauge condition (\ref{Colombo}) and boundary conditions  (\ref{gauge-bc})  are
\begin{eqnarray}
\label{Colombo-trans}& -\xi \theta_x + \eta= 0, \quad \theta(0)=\theta (L)= 0 .
\end{eqnarray}
Simplifying the equations in  (\ref{Stretch-E}) and the boundary conditions (\ref{a})-(\ref{aa}) by using (\ref{Colombo-trans}) yields
\begin{eqnarray}
\label{eqv}  &&\rho   \ddot v -{\alpha} v_{xx} -{\gamma}  \left((\phi^1)_x +  \dot\eta_x \right)= 0   \\
  \label{phi1} &&-\frac{{\e1} h^2}{12}(\phi^1)_{xx} +{\ep3} ~\phi^1 -{\gamma} v_x= \frac{\sigma_s (t) }{h}  \\
   \label{current-eq10} &&\frac{{\e1}h^2}{12}\ddot \theta +\mu \theta -\mu\eta_x+ \frac{{\e1}h^2}{12} (\dot \phi^1)_x=\frac{i_s^1(t)}{h} \\
\label{stretching} &&{\ep3}   \ddot \eta-\mu \eta_{xx} + \mu  \theta_x+\ep3 \dot \phi^1-{\gamma}\dot v_x=0
\end{eqnarray}
with the   boundary conditions
\begin{eqnarray}
\label{stretching-bcs}&&\left\{(\phi^1)_x (0)=\theta(0)=\eta_x (0)=\alpha  v_{x} (0) +{\gamma}\left(\phi^1+\dot\eta\right)\right\}_{x=0,L}=0.
\end{eqnarray}

It is shown in the next section that a well-posed system of equations has been obtained.

\section{Well-posedness}

\label{Sec-IV}

Consider first the existence and uniqueness of solutions to (\ref{stretching},\ref{stretching-bcs}) in the absence of control. It will be shown that these equations do have a unique solution. The solution defines a strongly continuous semigroup on a Hilbert space with norm corresponding to the physical energy. Moreover this semigroup is unitary, that is, the energy is conserved with time.

Define $\xi=\frac{{\e1}h^2}{12{\ep3}}.$ The elliptic equation (\ref{phi1}) with the associated boundary conditions can be written as
\begin{eqnarray}-\xi \phi^1_{xx} + \phi^1 = \frac{\gamma}{\ep3}v_x,\quad (\phi^1)_x(0)=(\phi^1)_x(L)=0. \label{elliptic}\end{eqnarray}
Consider
\begin{equation} -\xi D_x^2 \phi + \phi = z ,
\label{ellipticb}
\end{equation}
with $D_x^2 \phi = \phi_{xx}$ and domain
$${\rm Dom} (D_x^2) = \{ \phi \in H^2 (0,L) ,\quad \phi_x (0)=\phi_x (L)=0 \}. $$
Equation (\ref{elliptic}) has a unique solution for $\phi$ for any $z \in \cLtwo.$  Define the operator $P_\xi$
\begin{eqnarray}\label{Lgamma}P_{\xi}:=\left(-\xi D_x^2+I\right)^{-1}.\end{eqnarray}
It is well-known that ${P_\xi}$ is a compact operator on $\cLtwo$. Also,
 $P_\xi$ is a non-negative operator. To see this, let $P_\xi u=w.$ Then $w-\xi w_{xx}=u $ with $w_x(0)=w_x(L)=0,$ and
\begin{eqnarray}\label{non-neg}\left<P_\xi u, u\right>_{\cLtwo}=\left<w, w-\xi w_{xx}\right>_{\cLtwo}= \|w\|^2_{\cLtwo} + \xi \|w_{x}\|^2_{\cLtwo}\ge 0.
\end{eqnarray}

Therefore,  equation (\ref{phi1}) has the solution
\begin{eqnarray}\label{solutio}
\phi^1=\left\{
  \begin{array}{ll}
   \frac{{\gamma}}{{\ep3}}~P_{\xi} v_x , & \sigma_s(t)\equiv 0, i_s^1(t)\ne 0, \\
   \frac{{\gamma}}{{\ep3}}~P_{\xi} v_x + \frac{\sigma_s(t)}{\ep3 h} \left(H(x)-H(x-L)\right)+K , & \sigma_s(t)\ne 0, i_s^1(t)\equiv 0.
  \end{array}
\right.
 \end{eqnarray}
where $K$ is an arbitrary constant. In the case of  current actuation, the solution is unique.

Using (\ref{solutio}), the stretching equations (\ref{stretching}) are rewritten as
\begin{eqnarray}
& \rho  \ddot v-{\alpha} v_{xx} -\frac{{\gamma}^2}{{\ep3}}(P_{\xi} v_x)_x- {\gamma}\dot \eta_x=  \frac{\gamma\sigma_s(t)}{\ep3 h} \left(\delta(x)-\delta(x-L)\right)  & {\mbox{in}} ~\Omega \times \mathbb{R}^+ \\
 & \frac{\e1 h^2}{12} \ddot \theta +\mu  \theta -\mu \eta_x +\frac{{\e1}h^2}{12}\frac{{\gamma}}{{\ep3}}(P_{\xi} \dot v_x)_x=\frac{i_s^1(t)}{h}   & {\mbox{in}} ~\Omega \times \mathbb{R}^+\\
\label{IVP} &{\ep3}  \ddot \eta -\mu \eta_{xx} + \mu  \theta_x- {\gamma}\left(\dot v_x-(P_\xi \dot v_x)\right)= 0& {\mbox{in}}~\Omega \times \mathbb{R}^+~
\end{eqnarray}
with the same boundary conditions at $x=0, L$,
\begin{equation}
\label{BC}
{\alpha} v_{x} +\frac{{\gamma}^2 }{\ep3}{P_\xi} v_x + {\gamma} \dot \eta  =
\theta=\eta_{x}=0 .
\end{equation}
Note that the operator $P_\xi$ increases the mechanical stiffness in  the first equation of (\ref{IVP}). This stiffening does not occur if the potential $\phi$ in (\ref{scalarpot}) is assumed to vary linearly  with thickness, instead of quadratically as assumed here.

Defining the state variable,
$$\y=\begin{bmatrix} v_x \\ \theta \\ \eta \\ \dot v \\ \dot \theta \\ \dot \eta \end{bmatrix} = \begin{bmatrix} y_1 \\ y_2\\ \vdots \\y_6 \end{bmatrix}  $$
the natural energy associated with (\ref{IVP}) is
\begin{eqnarray}
  \nonumber\mathrm{E}(t)=\frac{1}{2}\int_0^L \left\{\rho  |y_4|^2 + \frac{\e1 h^2}{12} |y_5|^2 + \ep3 |y_6|^2 +  {\alpha} |y_1|^2 + \frac{{\gamma}^2  }{\ep3}({P_\xi} y_1) \bar y_1  \right.\\
\label{Energy-class} \left. + \mu  | y_2 -(y_3)_x|^2 \right\} dx, ~~~~t\in \mathbb{R}.
\end{eqnarray}
Writing
$$ H^1_0(0,L)=\{f\in H^1(0,L)~:~ f(0)=f(L)=0 \},$$
the energy  motivates definition of the
 linear space
    \begin{eqnarray}\nonumber
    \mathrm H &=  &\left\{ \y \in  \cLtwo \times H^1_0(0,L)\times {H}^1(0,L) \times  \cLtwo\times  \cLtwo \times   \cLtwo \right.,\\
\label{H-cur} && \left. -\xi(y_2)_x + y_3= -\xi (y_5)_x + y_6=0\right\}
    \end{eqnarray}
   and  bilinear form \begin{eqnarray}
\nonumber &&\left<\y, \z\right>_{\mathrm{H}}=\int_0^L \left\{\rho  y_4 \bar z_4  + \frac{\e1 h^2}{12} y_5\bar z_5 + \ep3  y_6\bar z_6 +  {\alpha} y_1 \bar z_1   + \frac{{\gamma}^2 }{\ep3}{P_\xi} y_1 \bar z_1 \right.\\
\label{inner-pro} &&\quad\quad\quad \left.   +\mu  y_2 \bar z_2
 + \mu  (y_3)_x (\bar z_3)_x -\mu  y_2  (\bar z_3)_x-\mu  ( y_3)_x \bar z_2 \right\} dx.~~~~\quad\quad
 \end{eqnarray}

\begin{thm}\label{complete}
The form (\ref{inner-pro}) defines an inner product on
the linear space $\mathrm H$. Moreover,  $E$ is the norm induced by this inner product and  $\mathrm H$ is complete.
\end{thm}

\textbf{Proof:}
It is straightforward to verify that (\ref{inner-pro}) defines a sesquilinear form.
The main problem is to show that this bilinear form  (\ref{inner-pro}) is coercive. This follows since $P_\xi$ is a self-adjoint positive operator on $\cLtwo.$
Using the gauge condition (\ref{Colombo-trans}), and Poinc\'{a}re's inequality with the Poinc\'{a}re constant $C,$
 \begin{eqnarray}\nonumber  - \int_0^L y_2 (\bar y_3)_x~dx  &=& -\xi\int_0^L y_2(\bar y_2)_{xx}~dx= \xi \int_0^L |(y_2)_{x}|^2 ~dx \ge C\xi \int_0^L |y_2|^2~dx
\end{eqnarray}
 Therefore, (\ref{inner-pro}) is a valid  inner product on $\mathrm H$ and so defines a norm.
 It is straightforward to verify that $E(t)$  as defined in (\ref{Energy-class}) is the norm induced by (\ref{inner-pro}).

It can also easily be shown that $\mathrm H$ with this norm is complete. This follows from the fact that the gauge constraints in $\mathrm H$ are satisfied weakly, i.e.
\begin{eqnarray}
\nonumber 0&=&\left<\xi y_2, \phi_x\right>_{\cLtwo}+ \left<\y_3,\phi\right>_{\cLtwo}\\
\nonumber 0&=&\left<\xi y_5, \phi_x\right>_{\cLtwo}+ \left<\y_6,\phi\right>_{\cLtwo}
\end{eqnarray}
for every $\phi\in H^1(0,L).$ Therefore a Cauchy sequence  $\{Y_n\}$ in $\mathrm H$ converges to $Y\in \mathrm H.$
  $\square$

Define  \begin{eqnarray}
\label{A1}&& A_1= \left( {\begin{array}{*{20}c}
 D_x \left(\frac{{\alpha}}{\rho} I+ \frac{{\gamma}^2}{\ep3 \rho} {P_\xi} \right) & 0 & 0    \\
    0 & \frac{-12\mu}{\e1 h^2}I & \frac{12\mu}{\e1 h^2} D_x \\
    0 & -\frac{\mu}{\ep3}  D_x & \frac{\mu}{\ep3}D_x^2  \\
\end{array}} \right)\\
\label{SSemigroup}&& A_2= \left( {\begin{array}{*{20}c}
   0 & 0 &  \frac{{\gamma}}{\rho}D_x   \\
    -\frac{{\gamma}}{ \ep3}D_x {P_\xi} D_x & 0 &0  \\
    \frac{{\gamma}}{\ep3}(I-{P_\xi}) D_x & 0 & 0   \\
\end{array}} \right), \quad \mc A= \left( {\begin{array}{*{20}c}
   0 & I_{3\times 3}  \\
   A_1 & A_2  \\
\end{array}} \right)\end{eqnarray}
with
\begin{eqnarray}\nonumber {\rm Dom}(\mc A)&=& \left(H^1(0,L) \times   H^2(0,L)\times H^2(0,L) \times H^1(0,L)\times H^1_0(0,L)\times (H^1(0,L))\right) \\
\label{domain-cur} && \bigcap\left\{ \y \in  \mathrm H  : \left|~\left(\frac{}{}{\alpha} I + \frac{{\gamma}^2}{\ep3} {P_\xi}\right) y_1 +{\gamma} y_6=y_2=(y_3)_x ~\right|_{x=0,L}=0 \right\} .
\end{eqnarray}

 \begin{Lem} \label{density}The operator $\mc A$ is densely defined in $ \mathrm H $.
 \end{Lem}

 \textbf{Proof:} Let $\y_n \in {\rm Dom}(\mc A)\to \y \in \mathrm H.$ Then
 \begin{eqnarray}\nonumber &&\|y_{1n}-y_1\|_{\cLtwo}\to 0 \\
\nonumber &&\|y_{2n}-y_2\|_{H^1_0(0,L)}\to 0 \\
\nonumber &&\|y_{3n}-y_3\|_{H^1(0,L)}\to 0 \\
\nonumber &&\|y_{4n}-y_4\|_{\cLtwo}\to 0 \\
\nonumber &&\|y_{5n}-y_5\|_{\cLtwo}\to 0 \\
\label{guzel}&& \|y_{6n}-y_6\|_{\cLtwo}\to 0.
\end{eqnarray}
For every $\z\in (C_0^\infty[0,L])', $ we have $\left<\y_n,z\right>\to \left<y, z\right>$ in $C_0^\infty[0,L].$ This also shows that the trace for $y_2$ is well defined.   For this we need the Green's formula as the following
For $y_{2n}\in H^1_0(0,L)$ and $\phi\in H^1(0,L)$
\begin{eqnarray}\left<{y_{2n}}_x,\phi\right>_{\cLtwo}&=& -\left<y_{2n},\phi_x\right>_{\cLtwo} + \left.y_{2n} \phi\right|_{x=0}^L \\
\nonumber &=& -\left<y_{2n},\phi_x\right>_{\cLtwo}\\
\nonumber &\to& -\left<y_{2},\phi_x\right>_{\cLtwo}\\
\nonumber &=& \left<y_{2x},\phi_x\right>_{\cLtwo} - \left.y_{2} \phi\right|_{x=0}^L\end{eqnarray}
and therefore $\left.y_2\right|_{x=0,L}=0.$ We check whether $\y$ satisfies the gauge conditions.
\begin{eqnarray}\nonumber 0=\left<-\xi (y_{2n})_x+ y_{3n},\phi\right>_{\cLtwo}&=&\left<\xi y_{2n},\phi_x\right>_{\cLtwo} + \left<y_{3n},\phi\right>_{\cLtwo}\\
\nonumber &\to& \left<\xi y_{2},\phi_x\right>_{\cLtwo} + \left<y_{3},\phi\right>_{\cLtwo}\\
\nonumber &=& \left<-\xi (y_{2})_x+y_{3},\phi\right>_{\cLtwo},
\end{eqnarray}
and,
\begin{eqnarray}\nonumber 0=\left<-\xi (y_{5n})_x+ y_{6n},\phi\right>_{\cLtwo} &\to& \left<\xi y_{5},\phi_x\right>_{\cLtwo} + \left<y_{6},\phi\right>_{\cLtwo}\\
\nonumber &=& \left<\xi y_{5},\phi_x\right>_{\cLtwo}+\left<y_{6},\phi\right>_{\cLtwo}
\end{eqnarray}
Now we show that  $\mc A : {\rm Dom} ({\mc A})  \to \mathrm H.$

\vspace{0.1in}
 \begin{Lem} \label{pxi} Let ${\rm Dom} (D_x^2)=\{w\in H^2(0,L) ~:~ w_x(0)=w_x(L)=0\}.$  The operator $\frac{1}{\xi}(P_\xi-I)$  is continuous, self-adjoint and non-positive on $\cLtwo.$ Moreover, for all $w\in {\rm Dom} (P_\xi), $
 $$ J=D_x^2 ~P_\xi= D_x^2 (I-\xi D_x^2)^{-1}w.$$
\end{Lem}
\vspace{0.1in}

\textbf{Proof:} 
 Define $J=\frac{1}{\xi}(P_\xi-I).$ Continuity and self-adjointness easily follow from the definition of $J.$ We first prove that $J$ is a non-positive operator. Let $u\in \cLtwo.$ Then $(I-\xi D_x^2)^{-1} u=s$ implies that $s\in {\rm Dom} (D_x^2)$ and $s-\xi s_{xx}=u$
\begin{eqnarray}\nonumber \left<Ju, u\right>_{\cLtwo}&=& \left<\frac{1}{\xi}(P_\xi-I) u, u\right>_{\cLtwo}\\
\nonumber &=& \frac{1}{\xi}\left< s - s + \xi s_{xx}, s-\xi s_{xx}\right>_{\cLtwo}\\
\nonumber &=& - \|  s_{x} \|^2_{\cLtwo} - \xi \|  s_{xx}\|^2_{\cLtwo}.
\end{eqnarray}
Let $Jw=\frac{1}{\xi}(P_\xi-I) w$ and $v:=P_\xi w.$ Then $v-\xi v_{xx}=w.$ By a simple rearrangement of the terms $$Jw=\frac{1}{\xi} (v-w)=\frac{1}{\xi}(v-v+\xi v_{xx})=v_{xx}=D_x^2 P_\xi w.~~ \square$$

\begin{Lem} \label{mapping} The operator $\mc A : {\rm Dom} {\mc A}  \to \mathrm H.$
\end{Lem}

\textbf{Proof:} Let $\y\in {\rm Dom}(\mc A).$ Then
$$\mc A \y=\left( \begin{array}{c}
y_4\\
y_5\\
y_6\\
\frac{{\alpha}}{\rho}(y_1)_{x} + \frac{{\gamma}^2}{\ep3 \rho} ({P_\xi} y_1)_x+\frac{{\gamma}}{\rho }( y_6)_x\\
-\frac{12\mu}{\e1 h^2} y_2 +\frac{12\mu}{ \e1 h^2} (y_3)_x-\frac{{\gamma}}{ \ep3}( {P_\xi} (y_4)_x)_x  \\
\frac{\mu}{\ep3}(y_3)_{xx}-\frac{\mu}{ \ep3} (y_2)_x +\frac{{\gamma}}{\ep3}(I-{P_\xi}) (y_4)_x
\end{array} \right). $$
First observe that since $\y \in {\rm Dom}(\mc A),$ we automatically have $y_4 \in{H}^1(0,L),$ $y_5 \in H^1_0(0,L),$ and $y_6 \in  {H}^1(0,L).$
Next,  $\frac{{\alpha}}{\rho}(y_1)_{x} + \frac{{\gamma}^2}{\ep3 \rho} ({P_\xi} y_1)_x-\frac{{\gamma}}{\rho }(y_3)_x \in  \cLtwo$ and $-\frac{12\mu}{\e1 h^2}y_2+\frac{12\mu}{\e1 h^2}(y_3)_x-\frac{{\gamma}}{ \ep3}({P_\xi} (y_4)_{x} )_x \in \cLtwo $
  follows from definition of $P_{\xi} .$  Finally, $\frac{\mu}{\ep3}(y_3)_{xx}-\frac{\mu}{ \ep3} (y_2)_x +\frac{{\gamma}}{\ep3}(I-{P_\xi})(y_4)_x \in  \cLtwo$
 also follows from  the definition of $P_\xi$.

  Since $\y \in {\mbox{Dom}}(\mc A),$ $-\xi(y_5)_x + y_6= 0.$ Next, we show that the other gauge constraint is satisfied. Using  Lemma \ref{pxi}, we obtain
\begin{eqnarray}&&-\frac{\e1 h^2}{12}\left[-\frac{12\mu}{\e1 h^2} y_2
+\frac{12\mu}{\e1 h^2}(y_3)_x-\frac{{\gamma}}{ \ep3}({P_\xi} (y_4)_{x})_x\right]_x\\
  \nonumber && \quad\quad\quad+\ep3 \left[\frac{\mu}{\ep3}(y_3)_{xx}-\frac{\mu}{ \ep3} (y_2)_x +\frac{{\gamma}}{\ep3}(I-{P_\xi}) (y_4)_x\right]\\
\nonumber &=& \mu  (y_2)_x-\mu  (y_3)_{xx}+\frac{\e1 {\gamma}h^2 }{12\ep3}({P_\xi} (y_4)_{x})_{xx} + \mu (y_3)_{xx} -\mu (y_2)_x  +{\gamma}(I-{P_\xi})(y_4)_x \\
\nonumber &=& \mu  (y_2)_x-\mu  (y_3)_{xx}+ \gamma \xi({P_\xi} (y_4)_{x})_{xx} + \mu (y_3)_{xx} -\mu (y_2)_x  +{\gamma}(I-{P_\xi})(y_4)_x \\
\nonumber &=&0. \square\end{eqnarray}

\vspace{0.1in}
 \begin{thm} \label{skew-adjoint-cur}The operator $\mc A$  satisfies $\mc{A}^*=-\mc{A}$ on  $\mathrm{H},$ and
  $\mathcal{A}: {\rm Dom} (\mc A) \subset \mathrm H \to \mathrm H$ defined by (\ref{SSemigroup})
is the generator of a unitary semigroup  $\{e^{\mathcal{A}t}\}_{t\ge 0}$ on $\mathrm H.$
 \end{thm}

\textbf{Proof:} Let $\y, \z \in {\rm Dom}(\mathcal{A}).$ Then we have
  \begin{eqnarray}
&& \nonumber \left<\mathcal{A}\y, \z\right>_{\mathrm{H}}  = \left<\left( \begin{array}{c}
y_4\\
y_5\\
y_6\\
\frac{{\alpha}}{\rho}(y_1)_{x} + \frac{{\gamma}^2}{\ep3 \rho} ({P_\xi} y_1)_x+\frac{{\gamma}}{\rho }( y_6)_x\\
-\frac{12\mu}{\e1 h^2} y_2+\frac{12\mu}{ \e1 h^2} (y_3)_x-\frac{{\gamma}}{ \ep3}( {P_\xi} (y_4)_x)_x  \\
\frac{\mu}{\ep3}(y_3)_{xx}-\frac{\mu}{ \ep3} (y_2)_x +\frac{{\gamma}}{\ep3}(I-{P_\xi}) (y_4)_x
\end{array} \right), \left( \begin{array}{c}
 z_1 \\
 z_2 \\
 z_3\\
 z_4\\
 z_5\\
 z_6\\
 \end{array} \right)\right>_{\mathrm{H}}\\
 \nonumber && = \int_0^L \left\{  \left({\alpha} (y_1)_{x} +
   \frac{{\gamma}^2 }{\ep3} ({P_\xi} (y_1)_x+ {\gamma} (y_6)_x\right) \bar z_4 +\left(-\mu  y_2 +\mu (y_3)_x-\gamma \xi( {P_\xi} (y_4)_x)_x
\right)\bar z_5\right. \\
    \nonumber &&~~ +  \left(\mu (y_3)_{xx}-\mu (y_2)_x +{\gamma} (I-{P_\xi}) (y_4)_x\right) \bar z_6
     + {\alpha}  y_4 \bar z_1  + \frac{{\gamma}^2}{\ep3} {P_\xi} y_4 \bar z_1  +  \mu  y_5 \bar z_2  + \mu  (y_6)_x (\bar z_3)_x \\
\nonumber && \left. ~~- \mu   y_5 (\bar z_3)_x - \mu  (y_6)_x \bar z_2 \right\}~dx\\
 \nonumber && = \int_0^L \left\{  -\left({\alpha} y_1 +
   \frac{{\gamma}^2 }{\ep3} {P_\xi} y_1+ {\gamma} y_6\right) (\bar z_4)_x +\left(-\mu  y_2 +\mu (y_3)_x-\gamma \xi( {P_\xi} (y_4)_x)_x
\right)\bar z_5\right.\\
    \nonumber &&~~ +  \left(\mu (y_3)_{xx}-\mu (y_2)_x +{\gamma} (I-{P_\xi}) (y_4)_x\right) \bar z_6
     + {\alpha}  (y_4)_x (\bar z_1)_x \\
\label{islem}  &&~~  \left. + \frac{{\gamma}^2}{\ep3} {P_\xi} (y_4)_x (\bar z_1)_x  +  \mu y_5 \bar z_2  + \mu  (y_6)_x (\bar z_3)_x - \mu  y_5(\bar z_3)_x - \mu  (y_6)_x \bar z_2 \right\}~dx
 \end{eqnarray}
where we integrated  the first integral by using the boundary condition (\ref{domain-cur}). Moreover, the integrals $-\gamma\xi \int_0^L ( {P_\xi} (y_4)_x)_x \bar z_5~ dx$  and $-\gamma \int_0^L P_\xi (y_4)_x \bar z_6~dx$ cancel each other by the gauge condition (\ref{H-cur}). Now we focus our attention to the last term in the first integral. We replace the term $ y_6(\bar z_4)_x $ by
$ y_6  P_\xi ^{-1}P_{\xi} (\bar z_4)_x$ since $P_{\xi}^{-1}P_{\xi} = (I-\xi D_x^2)P_\xi=I.$ By the gauge condition  (\ref{H-cur}) and integration by parts using (\ref{domain-cur})
\begin{eqnarray}\nonumber \gamma \int_0^L  y_6 (\bar z_4)_x~dx&=& \gamma\int_0^L \left( y_6 P_\xi (\bar z_4)_x
                 + y_6 (I-{P_\xi})(\bar z_4)_x \right)~dx\\
\nonumber &=& \gamma \int_0^L \left(\xi (y_5)_x  P_\xi (\bar z_4)_x + y_6  (I-{P_\xi})(\bar z_4)_x \right)~dx\\
\nonumber &=& \gamma\int_0^L \left( -\xi   y_5 (P_\xi (\bar z_4)_x)_x + y_6 (I-{P_\xi})(\bar z_4)_x \right)~dx.
\end{eqnarray}
 Integration by parts of the other terms in (\ref{islem}), using the gauge condition  (\ref{H-cur}) and  the boundary conditions (\ref{domain-cur}), results in
  \begin{eqnarray}
  \nonumber \left<\mathcal{A}\y, \z\right>_{\mathrm{H}}& =& \int_0^L -\left\{  \left({\alpha} \bar z_1 +
   \frac{{\gamma}^2 }{\ep3} {P_\xi} \bar z_1+ {\gamma} (\bar z_6)_x\right) y_4 \right. \\
  \nonumber &&~~+\left(-\mu  \bar z_2 +\mu (\bar z_3)_x-\gamma\xi ({P_\xi} (\bar z_4)_x)_x
\right) y_5\\
    \nonumber &&~~ +  \left(\mu (\bar z_3)_{xx}-\mu  (\bar z_2)_x +{\gamma} (I-{P_\xi}) (\bar z_4)_x\right)  y_6
     + {\alpha}  y_4 \bar z_1 \\
\nonumber &&\left.~~  + \frac{{\gamma}^2}{\ep3} {P_\xi} \bar z_4 y_1  +  \mu  (\bar z_5)(y_2)  + \mu  (\bar z_6)_x (y_3)_x
- \mu  (\bar z_5)(y_3)_x - \mu  (y_2) (\bar z_6)_x   \right\}dx.\\
\nonumber &=& \left<\y, -\mathcal{A} \z\right>_{{\mathrm H}}=\left<\y, \mathcal{A}^* \z\right>_{{\mathrm H}}.
\end{eqnarray}

 This implies that $\mc A$ is skew-symmetric.   To prove that $\mc A$ is skew-adjoint on $\mathrm H,$ i.e. $\mc A^*=-\mc A$ on $\mathrm H,$ with the same domains  it is required to show that  for any  $\g\in \mathrm{H}$  there is $\y \in  \text{Dom}(\mc A )$ so that   $\mc A \y=\g. $ This is equivalent to solving the  system of equations  for $ \y \in  \text{Dom}(\mc A ).$ Simplifying the equations leads to
  \begin{eqnarray}
\nonumber y_4=g_1\\
\nonumber  y_5=g_2\\
\nonumber  y_6=g_3\\
\nonumber \frac{{\alpha}}{\rho}(y_1)_{x} + \frac{{\gamma}^2}{\ep3 \rho} ({P_\xi} y_1)_x+\frac{{\gamma}}{\rho }( y_6)_x= g_4\\
\nonumber -\frac{12\mu}{\e1 h^2} y_2 +\frac{12\mu}{ \e1 h^2} (y_3)_x-\frac{{\gamma}}{ \ep3}( {P_\xi} (y_4)_x)_x =g_5\\
\label{den} \frac{\mu}{\ep3}(y_3)_{xx}-\frac{\mu}{ \ep3} (y_2)_x +\frac{{\gamma}}{\ep3}(I-{P_\xi}) (y_4)_x = g_6.
\end{eqnarray}

Since $\g\in  \mathrm H , $ we have
  $\g\in \cLtwo \times H^1_0(0,L)\times {H}^1(0,L) \times  \cLtwo\times  \cLtwo \times   \cLtwo$ and the components of $g$ satisfy the Gauge condition
  \begin{eqnarray}-\xi(g_2)_x + g_3=0, \quad  -\xi(g_5)_x + g_6=0.\end{eqnarray}
Then obviously, $-\xi(y_5)_x +  y_6=0.$ The solutions $y_2$ and $y_3$  are given by
\begin{eqnarray}
\nonumber  y_2(x)&=&(y_3)_x - \frac{\gamma\e1 h^2}{12\mu\ep3} ( {P_\xi} (g_1)_x)_x- \frac{\e1 h^2}{12 \mu }g_5\\
\label{sols} y_3(x)&=& P_\xi\left[-\frac{\gamma \e1 h^2}{12\mu\ep3}(P_\xi - I) (g_1)_x - \frac{\e1 h^2}{12 \mu} g_6\right].
\end{eqnarray}
 Now we check whether the gauge condition in (\ref{H-cur}) is satisfied
\begin{eqnarray}
\nonumber (y_2)_x&=& D_x^2(y_3) - \frac{\gamma\e1 h^2}{12\mu\ep3} D_x^2  {P_\xi} (g_1)_x - \frac{\e1 h^2}{12\mu }(g_5)_x\\
\nonumber &=& D_x^2 P_\xi\left[-\frac{\gamma \e1 h^2}{12\mu\ep3}(P_\xi -I) (g_1)_x - \frac{\e1 h^2}{12 \mu} g_6\right] - \frac{\gamma\e1 h^2}{12\mu\ep3} D_x^2  {P_\xi} (g_1)_x -\frac{\ep3}{\mu}g_6\\
\nonumber &=&\frac{1}{\xi}\left[ P_\xi -I\right]\left[-\frac{\gamma \e1 h^2}{12\mu\ep3}(P_\xi-I) (g_1)_x - \frac{\e1 h^2}{12 \mu} g_6\right] - \frac{\gamma\e1 h^2}{12\mu\ep3}\frac{1}{\xi}\left[P_\xi-I\right] (g_1)_x -\frac{\ep3}{\mu}g_6\\
\nonumber &=&\frac{1}{\xi}y_3
\end{eqnarray}
where we used (\ref{sols}) and Lemma \ref{pxi}. We find the unique solution $y_1\in \cLtwo$ by the Lax-Milgram theorem since
$$a(y_1,z_1)=\int_0^L \left[\frac{{\alpha}}{\rho} y_1  \bar z_1 + \frac{{\gamma}^2}{\ep3 \rho} {P_\xi} y_1 \bar z_1\right]~dx$$
is a continuous and coercive bilinear form on  $\cLtwo$  due to the positivity of the operator ${P_\xi}$ (see (\ref{non-neg})). Therefore $\y\in {\rm Dom}({\mc A}).$
This proves that $\mc A$ is skew-adjoint. By taking $\z=\y,$  (\ref{islem}) yields
  \begin{eqnarray}
\nonumber &&\left<\mc A \y, \y\right>_{ \mathrm H }= \int_0^L\left\{ {\alpha}  \left( y_4 {\bar y}_1 - {\bar y}_4  y_1 \right)
      +\frac{{\gamma}^2}{\ep3} \left({P_\xi} y_4 \bar y_1-{P_\xi} \bar y_4 y_1\right)\right. \\
      \nonumber && +{\gamma}\left((y_4)_x  \bar y_6 - (\bar y_4)_x y_6\right) + \mu  \left( y_5 \bar y_2- \bar y_5 y_2 \right)
      +\mu \left((y_3)_x \bar y_5 -(\bar y_3)_x  y_5  \right)\\
\nonumber  && \left. + \mu \left((\bar y_3)_x (y_6)_x-(y_3)_x (\bar y_6)_x\right) + \mu  \left( (\bar y_2)_x  y_6 - (y_2)_x \bar y_6\right) \right\}~dx.
\end{eqnarray}
Therefore ${\mbox{Re}}\left<\mc A \y, \y\right>_{ \mathrm H }={\mbox{Re}}\left<\mc A^* \y, \y\right>_{ \mathrm H }=0.$

It follows then from the Lumer-Phillips Theorem, that $\mc A$ generates a dissipative semigroup on ${\mathrm H}.$ Since $\mc A $ is skew-adjoint, the semigroup is unitary, that is in the absence of  control
$$\| y(t) \| = \| y(0)\|. \square$$

\vspace{0.1in}

The fact that $\mc A$ generates a unitary semigroup means that the norm and hence the energy $\mathrm{E}(t)$ is conserved along solution trajectories of (\ref{SSemigroup}) if there is no control term.

For $\y=[v_x,\theta,\eta,\dot v,\dot \theta,\dot \eta]^{\rm T}.$ The system (\ref{eqv})-(\ref{stretching-bcs}) can be written as
\begin{eqnarray}\label{system} \dot \y=\mc A\y + Bi_s(t), \quad \y(0)=y_0,
\end{eqnarray}
where the control operator $B$ is defined by $B\psi=[0~ 0~ 0~ 0~ \frac{12}{\varepsilon_1 h^3}~ 0]^{\rm T}=\psi_5.$
\begin{thm}\label{w-pf}
Let $T>0,$ and $i_s(t)\in \cLtwo(0,T).$ For any $\y_0 \in \mc H,$ $\y\in C[[0,T]; \mc H],$ and there exists a positive constants $c(T)$
such that (\ref{system}) satisfies
      \begin{eqnarray}\label{conc}\|\y (T) \|^2_{\mc H} &\le& c (T)\left\{\|\y^0\|^2_{\mc H} + \|i_s\|^2_{\cLtwo(0,T)}\right\}.
      \end{eqnarray}
\end{thm}
\vspace{0.1in}

\textbf{Proof:} The operator $\mc A: {\rm {Dom}}(\mc A) \to \mc H$ is a unitary semigroup by Lemma \ref{skew-adjoint-cur}. Therefore it is an  infinitesimal generator of $C_0-$semigroup of contractions by L\"umer-Phillips theorem. Since $$ \left< B i_s(t), \tilde \psi\right>_{\mathrm H,\mathrm H}=\int_0^L i_s(t)\psi_5(x) dx<\infty,$$ $B$ is an admissible control operator for the semigroup $\{e^{\mc A t}\}_{t\ge 0}$ corresponding to (\ref{system}), and hence the conclusion (\ref{conc})  follows. $\square$

\begin{thm}\label{iso}
The spectrum of $\mc A$ consists entirely of eigenvalues on the imaginary axis.
\end{thm}

\textbf{Proof:} By Lemma \ref{density}, ${\rm Dom}({\mc A})$ is densely defined  in ${\mathrm H}.$   ${\rm Dom}({\mc A})$ is also compact in ${\mathrm H}.$  To show this, let $\{Y_n\}$ be a bounded sequence in ${\rm Dom}({\mc A}),$ i.e.,
 $$\|y_{1n}\|_{H^1(0,L)},\|y_{2n}\|_{H^2(0,L)}, \|y_{3n}\|_{H^2(0,L)}, \|y_{4n}\|_{H^1(0,L)}, \|y_{5n}\|_{H^1(0,L)}, \|y_{6n}\|_{H^1(0,L)} <\infty.$$
From the Sobolev theory we know that both $H^1_0(0,L)$ and $ H^1(0,L)$ are compactly embedded in  $L^2(0,L).$ There exists a subsequence $\{Y_{n}\}\in \mathrm H,$  renamed similarly  as $\{Y_{n}\},$ such that
\begin{eqnarray}\nonumber &&\|y_{1n}-y_1\|_{\cLtwo}\to 0 \\
\nonumber &&\|y_{2n}-y_2\|_{H^1(0,L)}\to 0 \\
\nonumber &&\|y_{3n}-y_3\|_{H^1(0,L)}\to 0 \\
\nonumber &&\|y_{4n}-y_4\|_{\cLtwo}\to 0 \\
\nonumber &&\|y_{5n}-y_5\|_{\cLtwo}\to 0 \\
\label{dumb}&& \|y_{6n}-y_6\|_{\cLtwo}\to 0.
\end{eqnarray}
 Therefore for $\phi\in H^1(0,L)$
\begin{eqnarray}\nonumber 0=\left<-\xi (y_{2n})_x+ y_{3n},\phi\right>_{\cLtwo}&=&\left<\xi y_{2n},\phi_x\right>_{\cLtwo} + \left<y_{3n},\phi\right>_{\cLtwo}\\
\nonumber &\to& \left<\xi y_{2},\phi_x\right>_{\cLtwo} + \left<y_{3},\phi\right>_{\cLtwo}\\
\nonumber &=& \left<-\xi (y_{2})_x+y_{3},\phi\right>_{\cLtwo},
\end{eqnarray}
and,
\begin{eqnarray}\nonumber 0=\left<-\xi (y_{5n})_x+ y_{6n},\phi\right>_{\cLtwo} &\to& \left<\xi y_{5},\phi_x\right>_{\cLtwo} + \left<y_{6},\phi\right>_{\cLtwo}\\
\nonumber &=& \left<-\xi (y_{5})_x+y_{6},\phi\right>_{\cLtwo}
\end{eqnarray}
where we used (\ref{dumb}). This proves the compactness.  Moreover,  $0\in \rho({{\mc A}}) $ by Theorem \ref{skew-adjoint-cur}, it follows that $(\lambda I -{\mc A})^{-1}$
is compact at $\lambda=0,$ and thus compact for all $\lambda\in \rho({\mc A}).$ Hence the spectrum of ${\mc A}$ has all isolated eigenvalues. $\square$

\section{Conclusions}
In this paper a model for current actuation of a piezo-electric beam was derived in detail with fully dynamic electro-magnetic effects using Hamilton's Principle. An Euler-Bernoulli model was used for the mechanical model.  If  the Mindlin-Timoshenko small displacement assumptions are used instead, the bending equations  (\ref{Bending-E}) change substantially. However, stretching equations  in  (\ref{Stretch-E}) remain the same. Since the control only affects the stretching equations, the choice of beam model does not affect stabilizability.

With dynamic magnetic effects, the adjoint $B^*$ feedback in both the voltage- and current-controlled cases is electrical: for electrostatic models this feedback is mechanical.

 In the case of voltage actuation of a piezoelectric beam model  the control enters as a distribution.
Letting  $\delta$ indicate the Dirac delta function, $c$ a physical parameter,   the control operator is\cite{O-M}
$$B=c\left( \begin{array}{c}
 0_{3\times 1} \\
 \delta(x-L)-\delta(x)
 \end{array} \right)$$
 As when magnetic effects are neglected, the control operator $B$ is not bounded on the  state space $\mc H=(H^1(0,L))^2\times (\cLtwo)^2.$
However, when magnetic effects are included, the voltage-actuated piezo-electric beam is only exactly observable and exponentially stabilizable when the material parameters satisfy number-theoretical conditions. The system is asymptotically stabilizable under a wider set of parameter values \cite{O-M}. Explicit polynomial estimates for certain combination of parameters have been obtained  \cite{ozkan}.

Unlike voltage control,  for current actuation with magnetic effects,  the control operator  is bounded and rank 1. Thus,  it is compact and it is not possible to exponentially stabilize the piezoelectric beam;   see \cite{Gi80} or the textbook \cite{CZbook}.  Only asymptotic stabilization is possible.
 \begin{thm}  \label{main-thm2}
A given eigenvalue of $\mc A$  is asymptotically stabilizable if and only if the corresponding eigenfunctions  $\phi$  satisfy $\int_0^L \phi_5 dx \neq 0  .$
\end{thm}
\textbf{Proof:}
Let $\lambda $ be an eigenvalue of  $\tilde{\mc A}$ with eigenfunction $\psi,$ $ \| \psi \| = 1 :$
$$\lambda \psi = { \mc A }\psi - k_1 B B^* \psi. $$
If  $B^* \psi =0  $ then $\psi $ is an eigenfunction of $\mc A$ and so $B^* \psi \neq 0 . $ Then
$$\lambda = \langle { \mc A} \psi , \psi \rangle_H -k_1 | B^* \psi |_H^2$$
and
$$ \rm{Re} \, \lambda = -k_1  | B^* \psi |_H^2 < 0 .$$
Thus, since the spectrum consists only of imaginary eigenvalues (Theorem \ref{iso}),  and there are no eigenvalues on the imaginary axis,
Arendt-Batty's stability theorem \cite{A-B} implies  that the semigroup  is asymptotically stable. Conversely, if $\int_0^L \phi_5 dx =0$ for some eigenfunction, the corresponding
eigenvalue remains on the imaginary axis and the system is not asymptotically stable.
$\square$

With the state space $H$, based on energy, used in this paper, $0$ is an eigenvalue of $\mc A$ with an infinite-dimensional eigenspace
$$E=\{ y \in { H} ; \; y_4=y_5=y_6 =0 , \; y_3 \in  H^2 (0,L) , \;  y_{3x} (0)=y_{3x} (L)=0 , \;  y_2 = y_{3x}  \}  \subset {\mc D} ( \mc A ) .$$
Since $y_5=0 $ for all $y  \in E$, the $0$ eigenvalue is not stabilizable.
{This is typical for a structure with a rigid body mode; a voltage-controlled piezo-electric beam with natural boundary conditions, such as used here, also has an unstabilizable eigenvalue at $0.$}
 Further investigation is needed to determine whether the rest of the system is asymptotically stabilizable.

Stabilizability is quite different for electro-static models. For the voltage-controlled system,   an elliptic-type differential equation is obtained for charge, and once this equation is solved and back substituted to the mechanical equations, the system  reduces to a simple wave equation with the voltage control $V(t)$ acting at the free end of the beam.  This model is well-known to be exponentially stabilizable with  $B^*$ feedback, i.e. see \cite{Chen,ozkan3}.
 Hamilton's principle cannot be used to derive a  current-controlled system with electrostatic or quasi-static assumptions. Such a model can be obtained by adding   a circuit equation for the capacitance $\dot{V} = \frac{1}{C_p} i $ to the  voltage controlled model.
 The control operator is  bounded, so the system is not exponentially stabilizable. The $B^*$ feedback involves voltage.  The same analysis used for the voltage-controlled case in \cite{O-M} can be used to show that the system is  asymptotically stabilizable for certain parameter values.

 Charge actuation is mathematically very similar to voltage actuation for both electrostatic and quasi-static assumptions. This is because $\theta,\eta, \dot \theta,\dot \eta,\ll \phi^1.$ Without the terms $\theta,\eta, \dot \theta,\dot \eta$ in  (\ref{IVP})-(\ref{BC}), the model (\ref{IVP})-(\ref{BC}) for the clamped-free case
 becomes
 \begin{eqnarray}
 \label{electrostatic}  &&\left\{
  \begin{array}{ll}
 \rho  \ddot v-{\alpha} v_{xx} -\frac{{\gamma}^2}{{\ep3}}(P_{\xi} v_x)_x=\frac{\gamma\sigma_s(t)}{\ep3 h} \delta(x-L)  \quad {\mbox{in}} ~(0,L) \times \mathbb{R}^+ &\\
 \left. v(0)={\alpha} v_{x} +\frac{{\gamma}^2 }{\ep3}{P_\xi} v_x\right|_{x=L} =0, \quad t\in \mathbb{R}^+.&  \\
 v(x,0)=v_0(x), \quad \dot v(x,0)=v_1(x)\quad {\mbox{in}} ~(0,L)\end{array} \right.
\end{eqnarray}
Define $H^1_L(0,L)=\{\varphi\in H^1(0,L)~:~\varphi(0)=0\}.$
As for voltage control of the electrostatic  model, $B^*$ feedback control $\sigma_s(t)=-k B^*\y=-k \dot v(L,t)$ where $k>0,$ leads to an exponentially stable system; details can be found in \cite{ozkan2}.
The case of charge actuation with magnetic effects yields a model very similar to voltage control with magnetic effects and the system can be shown to be stabilizable for certain parameter values.

Thus, for all of voltage-, current and charge- control, magnetic effects have a significant effect on stabilizability.  Although the magnetic coupling $\mu$ is very small, stabilizability  of piezo-electric beams is qualitatively different for models with dynamic magnetics  than for an electrostatic models.


\begin{thebibliography}{9}
\bibitem{A-B} W. Arendt and C.J.K. Batty, {\sl{Tauberian theorems and stability of one-parameter semigroups,}} Trans. Amer. Math. Soc., 306 (1988), pp. 837--852.
\bibitem{Banks-Smith} H.T. Banks, R.C. Smith, Y. Wang, {\sl Smart material structures: Modelling, Estimation and Control}, Mason, Paris, 1996.

\bibitem{Cao-Chen} Y. Cao, X.B. Chen, {A survey of modeling and control Issues for piezo-electric actuators,} {\sl Journal of Dynamic Systems, Measurement, and Control, } {137-1} (2014), pp. 014001.

\bibitem{Chee} C.Y.K. Chee, L. Tong, and G.P. Steven, {\sl  A review on the modelling of piezoelectric sensors and actuators
incorporated in intelligent structures}, { J. Intell. Mater. Syst. Struct.}  {9} (1998), pp. 3–-19.
\bibitem{Chen} G. Chen, {\sl A note on the boundary stabilization of the wave equation,} { SIAM J. Control Optim.,}
19 (1981), pp. 106--113.
\bibitem{Comstock} R. Comstock, {\sl{Charge control of piezoelectric actuators to reduce hysteresis effects,}} (United States Patent \# 4,263,527, Assignee: The Charles Stark Draper Labrotary).
\bibitem{CZbook} R.F. Curtain and H.J. Zwart, {\sl{An Introduction to Infinite-Dimensional Linear Systems Theory,}}  Springer, Berlin, 1995.
	

\bibitem{L-D-I} R. Dautray and J.-L. Lions, {\sl{Mathematical Analysis and Numerical Methods for Science and Technology,}} Volume 1 Physical Origins and Classical Methods, Springer, Berlin, 1988.


\bibitem{Dest} Ph. Destuynder, I. Legrain, L. Castel, N. Richard, {\sl Theoretical, numerical and experimental discussion of the use of piezoelectric devices for control-structure interaction,} European J. Mech. A Solids, vol. 11 (1992), pp. 181--213.

 \bibitem{Duvaut-L} G. Duvaut, J.L. Lions, {\sl{Inequalities in Mechanics and Physics,}} Springer, Berlin, 1976.




\bibitem{Eringen} A.C. Eringen and G.A. Maugin, {\sl{Electrodynamics of Continua I, Foundations and Solid Media}}, Springer-Verlag, New York, 1990.
\bibitem{F-M} A.J. Fleming, S.O.R. Moheimani, {{\sl Precision current and charge
amplifiers for driving highly capacitive piezoelectric loads,}} Electron.
Lett., 39-3 (2003), pp. 282--284.

\bibitem{Gi80} J. S. Gibson, {\sl{A note on stabilization of infinite-dimensional linear oscillators by compact linear feedback}}, SIAM Journal of Control and Optimization, 18-3 (1980), pp. 311--316.

\bibitem{rob}  R.B. Gorbet and K.A. Morris and D.W.L. Wang,  {\sl{Passivity-based stability and control of hysteresis in smart actuators,}} {IEEE Transactions on Control Systems Technology}, 9-1 (2001),  pp. 5-16.
	
\bibitem{Hagood} N.W. Hagood, W.H. Chung, A.V. Flotov, {Modeling of piezoelectric actuator dynamics for active structural control,} J. Intelligent Material Systems and Structures, 1-3 (1990), pp. 327--354.



\bibitem{Hansen} S. W.  Hansen, {\sl Analysis of a plate with a localized piezoelectric patch}, Conference on Decision \& Control, Tampa, Florida, pp. 2952-2957 (1998).


   \bibitem{Lagnese-Lions} J.E. Lagnese, J.-L. Lions, {\sl Modeling Analysis and Control of Thin Plates,} Masson, Paris, 1988.




        \bibitem{Lee} P.C.Y. Lee, {\sl A variational principle for the equations of piezoelectromagnetism in elastic dielectric crystals,} Journal of Applied Physics, 69-11 (1991), pp. 7470--7473.



\bibitem{MGN} J.A. Main, E. Garcia and D.V. Newton, {{Precision position control of
piezoelectric actuators using charge feedback}}, {Journal of Guidance,
Control, and Dynamics,} 18-5 (1995), pp. 1068--73.


\bibitem{MG1} J.A. Main and E. Garcia, {\sl{Design impact of piezoelectric actuator
nonlinearities,}} { Journal of Guidance, Control, and Dynamics}, 20-2 (1997),
pp. 327--332.






\bibitem{O-M} K.A.  Morris, A.\"{O}. \"{O}zer, {\sl{Modeling and stabilizability of voltage-actuated piezoelectric beams with magnetic effects,}} SIAM Journal of Control and Optimization, 52-4 (2014), pp. 2371-2398.


    \bibitem{O-M3}  K.A. Morris, A.\"{O}. \"{O}zer, {\sl Comparison  of stabilization of  current-actuated  and voltage-actuated piezoelectric beams,}{ The $53^{\rm  rd}$ Proceedings of the IEEE Conf. on Decision \& Control,} Los Angeles, California, USA, pp. 571--576 (2014).


\bibitem{Newcomb-Flinn} C. Newcomb, I. Flinn, {{Improving the linearity of piezoelectric ceramic actuators,}} Electronic Letter, 18-11 (1982), pp. 442--443.
\bibitem{ozkan} A.\"O. \"Ozer, {\newblock{\sl Further stabilization and exact observability results for voltage-actuated piezoelectric beams with magnetic effects,}} Mathematics of Control, Signals, and Systems, 27-2 (2015), pp. 219--244.
\bibitem{ozkan2} A.\"O. \"Ozer, {\newblock{\sl Modeling and stabilization results for a charge or current-actuated active constrained layer (ACL) beam model with the electrostatic assumption,}} submitted to: Smart Structures NDE 2016: Active and Passive Smart Structures and Integrated Systems X, Piezo-based Materials and Systems II, Proc. of SPIE., arXiv: 1602.06368.


\bibitem{ozkan3} A.\"O. \"Ozer, {\newblock{\sl     Uniform stabilization of a multilayer Rao-Nakra sandwich beam,}} Evolution Equations and Control Theory, 4 (2013), pp. 695--710.



\bibitem{Rogacheva} N. Rogacheva, {\sl{ The Theory of Piezoelectric Shells and Plates}}, Boca Raton, FL: CRC Press, 1994.


\bibitem{Smith} R.C. Smith, {Smart Material Systems}, Society for Industrial and Applied Mathematics, 2005.

\bibitem{SmithII} R.C. Smith, C. Bouton and R. Zrostlik, {Partial and full inverse compensation for hysteresis
in smart material systems,} Proceedings of the American Control Conference, pp. 2750-2754, 2000.

\bibitem{Stanway} R. Stanway, J.A. Rongong, N.D. Sims, {Active constrained-layer damping: a state-of-the-art review,}  {\sl Automation \& Control Systems, }{217-6} (2003), pp. 437--456.


\bibitem{Tiersten} H.F. Tiersten, {\sl {Linear Piezoelectric Plate Vibrations} }, Plenum Press, New York, 1969.
\bibitem{Tzou} H.S. Tzou, {\sl{Piezoelectric shells, Solid Mechanics and Its applications 19,}} Kluwer Academic, The Netherlands, 1993.

   \bibitem{Yang1} J. Yang, {\sl{A review of a few topics in piezoelectricity,}} Appl. Mech. Rev., 59 (2006), pp. 335--345.

\end{thebibliography}
\end{document}